\documentclass[10pt,a4paper,twoside]{article}
\usepackage{amsfonts,amssymb,euscript,eufrak}
\usepackage{amscd}

\newcommand{\p}{\mathfrak{p}}
\newcommand{\tor}{\mathfrak{t}}\newcommand{\bo}{\mathfrak{b}}\newcommand{\un}{\mathfrak{u}}
\newcommand{\Q}{\mathbb{Q}}
\newcommand{\C}{\mathcal{C}}
\newcommand{\R}{\mathbb{R}}
\newcommand{\Z}{\mathbb{Z}}
\newcommand{\N}{\mathbb{N}}

\newcommand{\x}{\mathfrak{x}}
\newcommand{\al}{\alpha}

\newcommand{\ol}{\mathfrak{o}}

\newcommand{\db}{{\bf b}}
\newcommand{\dX}{{\mathfrak{X}}}
\newcommand{\Dh}{D(H,K)} \newcommand{\Dg}{D(G,K)}\newcommand{\Dgr}{D_r(G,K)}
\newcommand{\Do}{D(H_0,K)}\newcommand{\Dgo}{D(G_0,K)} \newcommand{\Dgor}{D_r(G_0,K)}

\newcommand{\Dor}{D_r(H_0,K)}

\newcommand{\gr}{gr\dot{}_r\,}\newcommand{\gor}{gr\dot{}\,}

\newcommand{\Co}{C^{an}(H_0,K)}

\newcommand{\lie}{\mathfrak{g}}
\newcommand{\hdh}{\mathfrak{h}}
\newcommand{\hd}{\lie^0}
\newcommand{\hs}{\lie_\sigma^0}

\newcommand{\nr}{\|.\|_{\bar{r}}}

\newcommand{\Rep}{\mathcal{R}}

\newcommand{\kap}{\kappa}

\newcommand{\Rpa}{{\rm Rep}_K^{a}(G)}

\newcommand{\Rpan}{{\rm Rep}_K^{a}(G_0)}
\newcommand{\Rp}{{\rm Rep}_K(G)}

\newcommand{\Rpn}{{\rm Rep}_K(G_0)}

\newcommand{\M}{\mathcal{M}}

\newcommand{\Can}{C^{an}}

\newcommand{\Dgs}{D(G_\sigma,K)}

\newcommand{\Dgos}{D(G_0',K)}\newcommand{\Dgss}{D(G',K)}

\newcommand{\Dp}{D(P,K)}\newcommand{\Dpr}{D_r(P,K)}\newcommand{\Dpo}{D(P_0,K)}\newcommand{\Dpor}{D_r(P_0,K)}

\newcommand{\Dd}{D}
\newcommand{\Ddr}{\Dd_r}

\newcommand{\Mgn}{{\mathcal {M}}(\Dgo)}
\newcommand{\Mg}{{\mathcal {M}}(\Dg)}
\newcommand{\Cgo}{\mathcal{C}_{G_0}}

\newcommand{\Bang}{{\rm Ban}_{G}^{adm}(K)}

\newcommand{\tr}{\varphi}\newcommand{\trr}{\varphi_r}

\newcommand{\Ti}{R^{i}F_{\Q_p}^L}

\newcommand{\Gl}{G_{L/\Q_p}}\newcommand{\Hl}{H_{L/\Q_p}}
\newcommand{\gl}{\lie_{L/\Q_p}}
\newenvironment{pr}{\it Proof:\rm}{\hfill $\Box$\newline\newline}
\newtheorem{theo}{Theorem}[section]
\newtheorem{prop}[theo]{Proposition}
\newtheorem{cor}[theo]{Corollary}

\newtheorem{lem}[theo]{Lemma}

\begin{document}
\title{Analytic vectors in continuous $p$-adic representations}
\author{Tobias Schmidt}
\maketitle
\begin{abstract} Given a compact $p$-adic Lie group $G$ over a
finite unramified extension $L/\Q_p$ let $\Gl$ be the product over
all Galois conjugates of $G$. We construct an exact and faithful
functor from admissible $G$-Banach space representations to
admissible locally $L$-analytic $\Gl$-representations that
coincides with passage to analytic vectors in case $L=\Q_p$. On
the other hand, we study the functor "passage to analytic vectors"
and its derived functors over general basefields. As an
application we compute the higher analytic vectors in certain
locally analytic induced representations.
\end{abstract}
\section{Introduction}
\footnote[0]{2000 Mathematics Subject Classification: 22E50, 11S99
(primary), 11F70 (secondary).

Keywords: $p$-adic Lie groups, representation theory, analytic
vector.}Recently, Schneider and Teitelbaum initiated a systematic
study of continuous representations of $p$-adic Lie groups into
$p$-adic topological vector spaces (cf. [ST1-6]). A central result
in this theory is that in case of a compact group $G$ over $\Q_p$
the algebra of locally analytic distributions on $G$ is a
faithfully flat extension of the algebra of continuous
distributions.  As a consequence, passage to analytic vectors
constitutes an exact and faithful functor $F_{\Q_p}$ from
admissible Banach space $G$-representations to admissible locally
analytic $G$-representations. Due to its properties $F_{\Q_p}$ is
a basic tool in a possible classification of admissible
topologically irreducible (unitary) Banach space representations
which is of particular interest in the realm of the $p$-adic
Langlands programme. It is therefore a natural question (raised by
J. Teitelbaum, cf. [T]) how to correctly generalize the above
results to groups over arbitrary base fields $\Q_p\subseteq L$.

Given a compact locally $L$-analytic group $G$ simple examples
show that the naive analogues of the above results do not hold
(for example, $F_L$ is not exact and often zero). The reason, as
we believe, is that the notion of a $K$-valued locally
$L$-analytic function depends on embedding the base field $L$ into
the coefficient field $K$. Consequently, we introduce, at least in
case $L/\Q_p$ Galois, the various restrictions of scalars
$G_\sigma$ of $G$ via $\sigma\in Gal(L/\Q_p)$ together with their
function spaces. Let $\Gl:=\prod_\sigma G_\sigma$. Denoting as
usual by $D^c(.,K)$ and $D(.,K)$ continuous and locally analytic
$K$-valued distributions respectively we construct a ring
extension
\[D^c(G,K)\rightarrow D(\Gl,K)\] which reduces to the former in case $L=\Q_p$ and is faithfully flat in case
$L/\Q_p$ is unramified (Thm. \ref{faith}). To obtain from this a
well-behaved generalization of $F_{\Q_p}$ we introduce for any
Banach space $G$-representation $V$ the subspace $V_{\sigma-an}$
of $\sigma$-analytic vectors whose formation is functorial in $V$.
Denoting by $\Bang$ and $\Rpa$ the abelian categories of
admissible Banach space and locally analytic representations of
$G$ over $K$ respectively we construct a functor
\[F: \Bang\rightarrow {\rm Rep}^{a}_K(\Gl)\] that enjoys, in case
$L/\Q_p$ unramified, the following properties (Thm. \ref{new}): it
is exact and faithful and coincides with $F_{\Q_p}$ in case
$L=\Q_p$. Given $V\in\Bang$ the representation $F(V)$ is strongly
admissible. Viewed as a $G_{\sigma^{-1}}$-representation $F(V)$
contains $V_{\sigma-an}$ as a closed subrepresentation and
functorial in $V$. The results obtained so far generalize two main
theorems of Schneider and Teitelbaum (cf. [S], Thm. 4.2/3) to
unramified extensions $L/\Q_p$.

The functor $F_L$ (over a general finite extension $L/\Q_p$) being
nevertheless an important construction we continue our work by
studying its derived functors. More generally, for an arbitrary
locally $L$-analytic group $G$ we study the functor "passage to
$L$-analytic vectors" $F^L_{\Q_p}$ from admissible locally
$\Q_p$-analytic representations to admissible locally $L$-analytic
representations. Then $F_L=F^L_{\Q_p}\circ F_{\Q_p}$ if $G$ is
compact and we may deduce left-exactness of $F_L$. It is unclear
at present whether the categories of admissible locally analytic
representations have enough injective objects. Nevertheless, we
prove that $F^L_{\Q_p}$ extends to a cohomological
$\delta$-functor $\Ti$ between admissible representations
vanishing in degrees $i>([L:\Q_p]-1)\,{\rm dim}_LG$. The functors
$\Ti$ turn out to be certain Ext-groups and satisfy
$R^{i}F_L=\Ti\circ F_{\Q_p}$ with $R^{i}F_L$ the right-derived
functors of $F_L$.

As an application we study the interaction of the $\delta$-functor
$\Ti$ with locally analytic induction. Let $P\subseteq G$ be a
closed subgroup (satisfying a mild extra condition). For all
$i\geq 0$ we obtain (Thm. \ref{interaction})
\[\Ti\circ {\rm Ind}_{P_0}^{G_0}={\rm
Ind}_P^G\circ \Ti\] as functors on finite dimensional locally
$\Q_p$-analytic $P$-representations. Here, $(.)_0$ refers to the
underlying $\Q_p$-analytic group. In case that $G$ equals the
$L$-points of a quasi-split connected reductive group over $L$ and
$P$ a parabolic subgroup we deduce from this an explicit formula
for the higher analytic vectors in principal series
representations of $G$.

~\\{\it Notations}. Let $|.|$ be the $p$-adic absolute value of
$\mathbb{C}_p$ normalized by $|p|=p^{-1}$. Let $\Q_p\subseteq
L\subseteq K\subseteq \mathbb{C}_p$ be complete intermediate
fields with respect to $|.|$ where $L/\Q_p$ is a finite extension
of degree $n$ and $K$ is discretely valued. Let $\ol\subseteq L$
be the valuation ring. $G$ always denotes a locally $L$-analytic
group with Lie algebra $\lie$. Their restriction of scalars to
$\Q_p$ are denoted by $G_0$ and $\lie_0$. For any field $F$ denote
by $Vec_F$ the category of $F$-vector spaces. For any ring $R$
denote by $\mathcal{M}(R)$ the category of right modules. Let
$\kap=1$ or $2$ if $p$ is odd or even respectively. We refer to
[NVA] for all notions from non-archimedean functional analysis.

~\\{\it Acknowledgements.} The author would like to thank Peter
Schneider for leading the author's attention to this problem. He
is also grateful to Matthew Emerton and Jan Kohlhaase for some
helpful remarks. Part of this work was done during a stay of the
author at the D\'epartement de Math\'ematiques, Universit\'e
Paris-Sud funded by the European Network "Arithmetic Algebraic
Geometry". The author is grateful for the support of both
institutions.

\section{Fr\'echet-Stein algebras}

In this section we recall and discuss two classes of
Fr\'echet-Stein algebras: distribution algebras and
hyperenveloping algebras. For a detailed account on abstract
Fr\'echet-Stein algebras as well as distribution algebras as their
first examples we refer to [ST5]. For all basic theory on uniform
pro-$p$ groups we refer to [DDMS]. Throughout this work all
indices $r$ are supposed to satisfy the technical conditions $r\in
p^\Q$,~$p^{-1}<r<1$ and $r\notin \{p^{\frac{-1}{p^h-p^{h-1}}},
h\in\N\}$.

A $K$-Fr\'echet algebra $A$ is called (two-sided) {\it
Fr\'echet-Stein} if there is a sequence $q_1\leq q_2 \leq...$ of
algebra norms on $A$ defining its Fr\'echet topology and such that
for all $m\in\N$ the completion $A_m$ of $A$ with respect to $q_m$
is a left and right noetherian $K$-Banach algebra and a flat left
and right $A_{m+1}$-module via the natural map $A_{m+1}\rightarrow
A_{m}$.

\begin{theo} {\bf (Schneider-Teitelbaum)}\label{ST1} Given a compact locally $L$-analytic group $G$ the
 algebra $\Dg$ of $K$-valued locally analytic distributions on $G$ is Fr\'echet-Stein.
\end{theo}

This is [ST5], Thm. 5.1. We recall the construction thereby fixing
some notation: choose a normal open subgroup $H_0\subseteq G_0$
which is a uniform pro-$p$ group. Choose a minimal set of ordered
generators $h_1,...,h_d$ for $H_0$. The bijective global chart
$\Z_p^d\rightarrow H_0$ for the manifold $H_0$ given by
\begin{equation}\label{order}(x_1,...,x_d)\mapsto
h_1^{x_1}...h_d^{x_d}\end{equation} induces a topological
isomorphism $\Co\simeq C^{an}(\Z_p^d,K)$ for the locally convex
spaces of $K$-valued locally analytic functions. In this
isomorphism the right-hand side is a space of classical Mahler
series and the dual isomorphism $D(H_0,K)\simeq D(\Z_p^d,K)$
therefore realizes $\Do$ as a space of noncommutative power
series. More precisely, putting
$b_i:=h_i-1\in\Z[G],~\db^\al:=b_1^{\al_1}...b_d^{\al_d}$ for
$\al\in\N_0^d$ the Fr\'echet space $\Do$ equals all convergent
series
\begin{equation}\label{ortho}
\lambda=\sum_{\al\in\N_0^d}d_\al\db^\al\end{equation} with
$d_\al\in K$ such that the set $\{|d_\al|r^{\kap |\al|}\}_\al$ is
bounded for all $r$. The family of norms $\|.\|_r$ defined via
$\|\lambda\|_r:=\sup_\al |d_\al|r^{\kap|\al|}$ defines the
Fr\'echet topology. They are multiplicative and the corresponding
completions $\Dor$ are $K$-Banach algebras exhibiting a
Fr\'echet-Stein structure on $\Do$. Choose representatives
$g_1,...,g_r$ for the cosets in $H\backslash G$ and define on
$\Dgo=\oplus_i\, \Do\,g_i$ the norms
$||\sum_i\lambda_i\,g_i||_r:=\max_i\,||\lambda_i||_r$. The
completions $\Dgor$ are the desired Banach algebras for $\Dgo$.
Finally, $\Dg$ is equipped with the corresponding quotient norms
$\nr$ coming from the quotient map $\iota': \Dgo\rightarrow\Dg$.
The latter arises as the dual map to the embedding
\begin{equation}\label{emb}\iota: \Can(G,K)\subseteq\Can(G_0,K).\end{equation}
Passing to the norm
completions $\Dgr$ yields the appropriate Banach algebras.

We mention another important feature of $\Dg$ in case $G_0$ is a
uniform pro-$p$ group. Each algebra $\Dgor$ carries the filtration
defined by the additive subgroups\[F^s_r\Dgor:=
\{\lambda\in\Dgor,~||\lambda||_r\leq p^{-s}\}\] and
$F^{s+}_r\Dgor$ (defined as $F^s_r\Dgor$ via replacing $\leq$ by
$<$) for $s\in\R$. Put
\[\gr\Dgor:=\oplus_{s\in\R}\,F^s_r\Dgor/F^{s+}_r\Dgor\]
for the associated graded ring. Given $\lambda\in
F^s_r\Dgor\setminus F^{s+}_r\Dgor$ we denote by
$\sigma(\lambda)=\lambda + F_r^{s+}\Dgor\in\gr\Dgor$ the {\it
principal symbol} of $\lambda$. Note that $\Dgr\simeq\Dgor/I_r$ is
endowed with the quotient filtration where $I:={\rm
ker}~\Dgo\rightarrow\Dg$ and $I_r$ denotes the closure of
$I\subseteq\Dgor$. These filtrations are exhaustive, separated,
complete and quasi-integral (in the sense of [ST5], \S1). For
$\lambda\neq 0$ in $\Dgr$ the principal symbol
$\sigma(\lambda)\in\gr\Dgr$ is defined analogously.
\begin{theo}{\bf (Schneider-Teitelbaum)}\label{ST2}
If $G_0$ is a $d$-dimensional locally $\Q_p$-analytic group and
uniform pro-$p$ then there is an isomorphism of $\gor K$-algebras
\[\gr\Dgor\stackrel{\sim}{\longrightarrow} (\gor K)[X_1,...,X_d],~\sigma(b_i)\mapsto
X_i.\]
\end{theo}
This is [loc.cit.], Thm. 4.5. Since $\gr\Dgr$ equals a quotient of
$\gr\Dgor$ each $D_r(G,K)$ is a complete filtered ring with
noetherian graded ring and hence, is a Zariski ring (cf. [LVO],
II.2.2.1).

Working over the base field $L$ we shall need to impose a mild
additional condition on uniform subgroups. Let $G$ be a
$d$-dimensional locally $L$-analytic group whose underlying
$\Q_p$-analytic group is uniform. Any minimal ordered set of
generators for $G_0$ defines a global chart $\Z_p^{nd}\rightarrow
G_0$ and hence, determines a $\Q_p$-basis of $\lie_0$. Note that
$\lie_0\simeq\lie_L$ canonically over $\Q_p$ ([B-VAR], 5.14.5). We
call $G$ {\it uniform*} if the generators can be chosen in such a
way that this basis has the form $v_i\x_j$ for a $\Z_p$-basis
$v_1=1,v_2,...,v_n$ of $\ol_L$ and an $L$-basis $\x_1,...,\x_d$ of
$\lie_L$. By [Sch], Cor. 4.4 each locally $L$-analytic group has a
fundamental system of neighbourhoods of the identity consisting of
normal uniform* subgroups (note that being uniform* implies the
condition (L) used in [loc.cit.])

We turn to another closely related class of Fr\'echet-Stein
algebras. Let $G$ be a locally $L$-analytic group of dimension
$d$. Let $U(\lie)$ be the enveloping algebra and let $\Can_1(G,K)$
be the stalk at $1\in G$ of the sheaf of $K$-valued locally
$L$-analytic functions on $G$. It is a topological algebra with
augmentation whose underlying locally convex $K$-vector space is
of compact type (for the basic properties of such spaces we refer
to [ST2], \S1). Denote by $U(\lie,K):=\Can_1(G,K)'_b$ its strong
dual, the {\it hyperenveloping algebra} (cf. [P], \S8). The
notation reflects that, up to isomorphism, $U(\lie,K)$ depends
only on $\lie$. It is a topological algebra with augmentation on a
nuclear Fr\'echet space. There is a canonical algebra embedding
$U(\lie)\subseteq U(\lie,K)$ with dense image compatible with the
augmentations. The formation of $\Can_1(G,K)$ and $U(\lie,K)$ (as
locally convex topological algebras) is functorial in $G$ and
converts direct products into (projectively) completed tensor
products taken over $K$. Dualizing the strict surjection
$\Can(G,K)\rightarrow \Can_1(G,K)$ yields an injective continuous
algebra map $U(\lie,K)\rightarrow D(G,K)$ which is a topological
embedding with closed image.
\begin{theo}\label{jan} {\bf (Kohlhaase)} Suppose $G_0$ is a uniform pro-$p$ group.
Denote by $h_{11},...,h_{nd}$ a minimal set of ordered generators
and put $b_{ij}=h_{ij}-1$. Denote by $U_r(\lie,K)$ the closure of
$U(\lie,K)\subseteq\Dgr$. There is a number $\epsilon(r,p)\in\N$
depending only on $r$ and $p$ such that the (left or right)
$U_r(\lie,K)$-module $\Dgr$ is finite free on the basis
$\Rep:=\{\db^\al,~ \al_{ij}<\epsilon(r,p)
{\rm~for~all~}(i,j)\in\{1,...,n\}\times\{1,...,d\}\}$. Letting
$\x_1,...,\x_d$ be an $L$-basis of $\lie$ and
$\dX^\beta:=\x_1^{\beta_1}\cdot\cdot\cdot\x_d^{\beta_d}$ one has
as $K$-vector spaces
\[U_r(\lie,K)=\{\sum_{\beta\in\N_0^d}d_\beta\dX^\beta,~d_\beta\in K,~||d_\beta\dX^\beta||_{\bar{r}}\rightarrow 0
{\rm~for~}|\beta|\rightarrow\infty\}\] where the power series
expansions are uniquely determined. The Banach algebras
$U_r(\lie,K)$ exhibit $U(\lie,K)$ as a Fr\'echet-Stein algebra.
\end{theo}
\begin{pr}
This is extracted from (the proof of) [Ko1], Thm. 1.4.2. The
number $\epsilon(r,p)$ equals the (unique) value of $t$ where the
supremum $\sup_{t\in\N} |1/t|r^{\kap t}$ is attained.
\end{pr}
In [loc.cit.] the noetherian and the flatness property of the
family $U_r(\lie,K)$ is immediately deduced from the commutative
diagram \[
\begin{CD}
U_r(\lie,K)@>>>
U_{r'}(\lie,K)\\
@VVV    @VVV\\
D_r(G,K)@>>>D_{r'}(G,K)
\end{CD}
\]
for $r'\leq r$ in which the lower horizontal arrow is a flat map
between noetherian rings and the vertical arrows are, by the first
statement in the theorem, finite free ring extensions. We also
remark that the first statement in the theorem in case $L=\Q_p$ is
due to H. Frommer ([F], 1.4 Lem. 3, Cor. 1/2/3).
\begin{prop}\label{topgen0}\label{tobe}
Suppose $G_0$ is uniform* and let $\x_1,...,\x_d$ be the
corresponding basis of $\lie_L$. Endowing the extension
$U_r(\lie,K)\subseteq D_r(G,K)$ with the $\nr$-norm filtration the
map $\gr U_r(\lie,K)\rightarrow\gr\Dgr$ is finite free on the
basis $\sigma(\Rep)$. Moreover, $\gr U_r(\lie,K)$ equals a
polynomial ring in $\sigma(\x_1),...,\sigma(\x_d)$ and
$||.||_{\bar{r}}$ is multiplicative on $U_r(\lie,K)$. For any
$\lambda=\sum_{\beta\in\N_0^d}d_\beta\dX^\beta\in U_r(\lie,K)$ one
has
\[||\lambda||_{\bar{r}}=\sup_\beta
|d_\beta|c_r^{|\beta|}\] where $c_r\in\R_{>0}$ depends only on $r$
and $p$.
\end{prop}
\begin{pr}
Let $v_1,...,v_n$ be a corresponding $\Z_p$-basis of $\ol_L$ for
the uniform* group $G$. Then $h_{ij}:=\exp(v_i\x_j)$ are a minimal
set of topological generators for $G$. Put as usual
$b_{ij}:=h_{ij}-1$. Now $v_i\x_j=\log(1+b_{ij})$ is a $\Q_p$-basis
for $\lie_0$ and a short calculation yields
$\sigma(v_i\x_j)=\sigma(b_{ij})^{p^h}$ with $h$ depending only on
$r$ and $p$. Hence, Thm. \ref{ST2} translates the map $\gr
U_r(\lie_0,K)\rightarrow\gr\Dgor$ into the inclusion $(\gor
K)[X^{p^h}_{11},...,X^{p^h}_{nd}]\subseteq (\gor
K)[X_{11},...,X_{nd}]$. For $L=\Q_p$ we obtain from this all
statements together with the fact that the $U_r(\lie_0,K)$-module
basis $\Rep$ for $\Dgor$ is in fact orthogonal with respect to
$||.||_r$. By [Sch], Lem. 5.3/Prop. 5.5 the graded ideal $\gr I_r$
where $I_r={\rm ker}~(\Dgor\rightarrow\Dgr)$ is generated by the
elements $X_{ij}^{p^h}-\bar{v}_iX_{1j}^{p^h}$ where
$\bar{v}_i\in\gor K$ equals the residue class of $v_i$. By similar
arguments the same holds true for $\gr J_r$ where $J_r={\rm
ker}~(U_r(\lie_0,K)\rightarrow U_r(\lie,K))$. It follows that
$I_r=\oplus_{g\in\Rep}\,J_r\,g$. By orthogonality the quotient
norm on $U_r(\lie,K)$ with respect to $||.||_r$ and
$U_r(\lie_0,K)\rightarrow U_r(\lie,K)$ equals precisely $\nr$. In
other words, $\gr U_r(\lie_0,K)/gr J_r\simeq\gr U_r(\lie,K)$ and
since $\gr I_r\cap \gr U_r(\lie_0,K)=\gr J_r$ the first statement
follows. Now $\gr U_r(\lie_0,K)/gr J_r$ is readily seen to be a
polynomial ring in the residue classes of the
$X^{p^h}_{1j},~j=1,...,d$ which correspond to $\sigma(\x_j)\in\gr
U_r(\lie,K)$. This implies that $\nr$ is multiplicative on
$U_r(\lie,K)$, that the topological $K$-basis
$\dX^\beta,~\beta\in\N_0^d$ for $U_r(\lie,K)$ is in fact an
orthogonal basis with respect to $\nr$ and that
$c_r:=||\x_j||_{\bar{r}}=||\x_j||_r=||\log(1+b_{1j})||_r=\sup_{t\in\N}|1/t|r^{\kap
t}$ depends only on $r$ and $p$.\end{pr}

%%%%%%%%%%%%%%%%%%%%%%%%%%%%%%%%%%%%%%%%%%%%%%%%%%%%%%%%

Next we prove a proposition on the compatibility of two
Fr\'echet-Stein structures. This will be used in the proof of Thm.
\ref{interaction}.
\begin{lem}\label{topgen}
Let $G$ be a compact locally $L$-analytic group of dimension $d$
and $P\subseteq G$ a closed subgroup of dimension $l\leq d$. There
is an open normal subgroup $G'\subseteq G$ with the following
properties: it is uniform* with respect to bases $\x_1,...,\x_d$
and $v_1,...,v_n$. Furthermore, $P':=P\cap G'$ is uniform* with
respect to the bases $\x_1,...,\x_l$ and $v_1,...,v_n$.
\end{lem}
\begin{pr}
Denote the Lie algebras of $G$ resp. $P$ by $\lie_L$ resp.
$\mathfrak{p}_L$. Denote by $G_0, P_0$ the underlying locally
$\Q_p$-analytic groups. Applying [DDMS], Prop. 3.9, Thm. 4.2/4.5
we see that $P$ contains a uniform subgroup $P_1$ such that every
open normal subgroup of $P$ lying in $P_1$ is uniform itself.
After this preliminary remark we choose, according to [Sch], Cor.
4.4 and (the proof of) [loc.cit.], Prop. 4.3, a locally
$L$-analytic group $G'$ open normal in $G$ with the following
properties: it is uniform* with an $L$-basis $\x_1,...,\x_d$ of
$\lie_L$ such that $\x_1,...,\x_l$ is an $L$-basis of
$\mathfrak{p}_L$. Furthermore, we may arrange that $P':=G'\cap
P\subseteq P_1$. Since $\x_1,...,\x_l$ is an $L$-basis of
$\mathfrak{p}_L$ and $\exp$ may be viewed an exponential map for
$P$ the $nl$ elements $\exp(v_i\x_j),~i=1,...,n,~j=1,...,l$ are
part of a minimal generating system for $G'$ and lie in $G'\cap
P=P'$. Since they are pairwise different modulo $G'^p$, hence
modulo $P'^p$ it follows from ${\rm dim}_LP'=l$ that they form a
minimal generating system for the uniform group $P'$. Thus, $P'$
is uniform* with the required bases.
\end{pr}
\begin{prop}\label{family}
Let $G$ be a compact locally $L$-analytic group. There is a family
of norms $(||.||_r)$ on $\Dgo$ with the following properties: it
defines the Fr\'echet-Stein structure on $\Dgo$ as well as on the
subalgebra $\Dpo$. For each $r$ the completion $\Dgor$ is flat as
a $\Dpor$-module. The family of quotient norms $(\nr)$ defines the
Fr\'echet-Stein structure on $\Dg$ as well as on the subalgebra
$\Dp$. For each $r$ the completion $\Dgr$ is flat as a
$\Dpr$-module.
\end{prop}
\begin{pr}
Apply the preceding lemma to $P\subseteq G$ to find an open normal
subgroup $G'\subseteq G$ which is uniform* with respect to bases
$\x_j$ and $v_i$. The $nd$ elements $\exp(v_i\x_j)$ are then
topological generators for $G'$ where the first $nl$ elements
($l:={\rm dim}_LP$) generate the uniform group $P':=P\cap G'$.
Endow $D(G_0,K)$ and $D(G,K)$ with the Fr\'echet-Stein structures
constructed in the beginning of this section and restrict the
norms to $D(P_0,K)$ and $D(P,K)$. Then [ST6], Prop. 6.2. yields
all statements over $\Q_p$. By definition, the quotient map
$\Dgo\rightarrow\Dg$ restricts to the quotient map
$\Dpo\rightarrow\Dp$ whence it is easy to see that the restricted
norms $\nr$ on $\Dp$ equal the quotient norms. Hence they realize
a Fr\'echet-Stein structure on $\Dp$ and only the last claim
remains to be justified. By the argument given at the end of
[loc.cit.] it suffices to prove it for the pair $P'\subseteq G'$.
Applying Prop. \ref{tobe} to $D(P',K)$ and $D(G',K)$ we obtain a
commutative diagram of commutative algebras
\[
\begin{CD}
\gr U_r(\p,K)@>>>
\gr U_{r}(\lie,K)\\
@VVV    @VVV\\
\gr D_r(P',K)@>>>\gr D_{r}(G',K)
\end{CD}
\]
in which the vertical arrows are finite free ring extensions on
bases $\sigma (\Rep(\p))$ resp. $\sigma(\Rep (\lie))$. By our
assumptions and by the explicit shape of these bases there is a
set $S\subseteq\sigma(\Rep (\lie))$ such that $\sigma(\Rep
(\lie))=\{s\,t, s\in S,~t\in \sigma (\Rep(\p))\}$. It follows that
the map of $\gr D_r(P',K)$-modules
\[
\oplus_{g\in S}~\gr D_r(P',K)\otimes_{\gr U_r(\p,K)}\gr
U_r(\lie,K)\rightarrow\gr D_{r}(G',K)\] induced by
$(\lambda\otimes\mu)_g\mapsto\sum_g \lambda\,\mu\, g$ is
bijective. Using [ST5], Prop. 1.2 we are hence reduced to prove
the flatness of the upper horizontal arrow. But this equals the
inclusion of a polynomial ring over $\gor K$ in $\l$ variables
into one of $d$ variables (again by Prop. \ref{tobe}) which is
clearly flat.
\end{pr}
In case $L=\Q_p$ this result is precisely [ST6], Prop. 6.2. The
proof of our proposition was simplified by a remark of J.
Kohlhaase.

We finish this section with some results on the Lie algebra
cohomology of $U(\lie,K)$. Recall the homological standard complex
of free $U(\lie)$-modules $U(\lie)\otimes_L\dot{\bigwedge}\lie$
whose differential is given via \[
\begin{array}{r}
  \partial(\lambda\otimes\x_1\wedge...\wedge\x_q) =\sum_{s<t}(-1)^{s+t}\lambda\otimes
[\x_s,\x_t]\wedge\x_1\wedge...\wedge
\widehat{\x_s}\wedge...\wedge\widehat{\x_t}\wedge...\wedge\x_q  \\
     \\
  +\sum_{s}(-1)^{s+1}\lambda\x_s\otimes\x_1\wedge...\wedge
\widehat{\x_s}\wedge...\wedge\x_q.\\
\end{array}
\]
Composing with the augmentation $U(\lie)\rightarrow L$ yields a
finite free resolution of the $U(\lie)$-module $L$ and
$H^*(\lie,V):=h^*({\rm Hom}_{L}(\dot{\bigwedge}\lie,V))$ resp.
$H_*(\lie,V):=h_*(V\otimes_L\dot{\bigwedge}\lie)$ as objects in
$Vec_L$ for any $\lie$-module $V$. Assume $V\in Vec_K$ is nuclear
Fr\'echet or of compact type such that $\lie$ acts by continuous
$K$-linear operators. Endow each $V\otimes_L\bigwedge^q\lie$ resp.
${\rm Hom}_L(\bigwedge^q\lie,V)$ with the projective tensor
product topology resp. the strong topology. The obvious map
$V'_b\otimes_L\bigwedge^q\lie\rightarrow ({\rm
Hom}_L(\bigwedge^q\lie,V))'_b$ is a topological isomorphism and
identifies $V'_b\otimes_L\dot{\bigwedge}\lie$ with the strong dual
of ${\rm Hom}_L(\dot{\bigwedge}\lie,V)$ (e.g. [P], 1.4). Endow
$H^*(\lie,V)$ resp. $H_*(\lie,V)$ always with the induced
topologies.
\begin{lem}\label{isodual}
Let $V$ be a nuclear Fr\'echet space with continuous
$\lie$-action. Suppose the differential in
$V\otimes_L\dot{\bigwedge}\lie$ is strict. There are isomorphisms
of locally convex $K$-vector spaces
\begin{equation}
 H^*(\lie,V'_b)\simeq H_*(\lie,V)'_b\end{equation} natural in
$V$.
\end{lem}
\begin{pr}
Since $V\otimes_L\dot{\bigwedge}\lie$ consists of Fr\'echet spaces
the differential has closed image. By [ST2], Thm. 1.1, Prop. 1.2
the complex $(V\otimes_L\dot{\bigwedge}\lie)'_b$ consists of
spaces of compact type and has a strict differential with closed
image. Thus, we may substitute in the proof of [Ko2], Lem. 3.6 all
weak topologies by the strong topologies and obtain $K$-linear
bijections
$H_*(\lie,V)'_b=(h_*(V\otimes_L\dot{\bigwedge}\lie))'_b\simeq
h^*((V\otimes_L\dot{\bigwedge}\lie)'_b)$ which are readily seen to
be topological. Since ${\rm Hom}_{L}(\dot{\bigwedge}\lie,V'_b)$
consists of spaces of compact type the remark preceding the lemma
implies $(V\otimes_L\dot{\bigwedge}\lie)'_b\simeq {\rm
Hom}_{L}(\dot{\bigwedge}\lie,V'_b)$ topologically whence the
claim.
\end{pr}
\begin{prop}\label{acyclic}
One has $\oplus_*H_*(\lie,U(\lie,K))=H_0(\lie,U(\lie,K))=K$.
\end{prop}
\begin{pr} Over the complex numbers this follows from [P], Thm. 8.6.
In our setting our results allow to give a proof along the lines
of [ST6], Prop. 3.1. The case $*=0$ is clear. Now being a
Fr\'echet-Stein algebra (Thm.\ref{jan}) the topology on
$U(\lie,K)$ is nuclear Fr\'echet and the differential in
$U(\lie,K)\otimes_L\dot{\bigwedge}\lie$ is strict ([ST5], \S3). By
Lem. \ref{isodual} it suffices to prove
$H^*(\lie,C^{an}_1(G,K))=0$ for $*>0$. By [BW], VII.1.1. the
complex ${\rm Hom}_{L}(\dot{\bigwedge}\lie,C^{an}_1(G,K))$ equals
(up to sign) the stalk at $1\in G$ of the deRham complex of
$K$-valued global locally-analytic differential forms on the
manifold $G$. By the usual Poincar\'e lemma the latter is acyclic.
\end{pr}

\section{Continuous representations and analytic vectors}

We recall some definitions and results from continuous
representation theory relying on [S]. We introduce the notion of
analytic vector and prove some basic properties.

A {\it locally analytic $G$-representation} is a barrelled locally
convex Hausdorff $K$-vectorspace $V$ equipped with a $G$-action
via continuous operators such that for all $v\in V$ the orbit map
$o_v: G\rightarrow V,~g\mapsto g^{-1}v$ lies in $\Can(G,V)$, the
space of $V$-valued locally analytic functions on $G$. With
continuous $K$-linear $G$-maps these representations form a
category $\Rp$. Endowing $\Can(G,V)$ with the left regular action
($(g.f)(h)=f(g^{-1}h)$) yields $\Can(G,V)\in\Rp$ and a
$G$-equivariant embedding $o: V\rightarrow\Can(G,V),~v\mapsto
o_v.$ Let $\Rpa\subseteq\Rp$ be the full subcategory of {\it
admissible} representations. Denote the abelian category of
coadmissible modules by $\C_G$. There is an anti-equivalence
$\Rpa\simeq \C_G$ via $V\mapsto V'_b$. In particular, any
$M\in\C_G$ has a nuclear Fr\'echet topology (the {\it canonical}
topology). If $G$ is compact $\C_G$ contains all finitely
presented modules. In this case, a $V\in\Rpa$ such that $V'$ is
finitely generated is called {\it strongly admissible}.

Now let $G$ be compact and $K/\Q_p$ be finite. A {\it Banach space
representation} of $G$ is a $K$-Banach space $V$ with a linear
action of $G$ such that $G\times V\rightarrow V$ is continuous.
Let $D^c(G,K)$ be the algebra of continuous $K$-valued
distributions on $G$. Denote by $\Bang$ the abelian category of
admissible representations and by $\mathcal{M}^{fg}(D^c(G,K))$ the
finitely generated modules. There is an anti-equivalence
$\Bang\simeq\mathcal{M}^{fg}(D^c(G,K))$ via $V\mapsto V'$. In
particular, $\Bang$ has enough injective objects. If $V\in\Bang$
then $v\in V$ is called a {\it locally $L$-analytic} vector if the
orbit map $g\mapsto gv$ lies in $\Can(G,V)$. The subspace
$V_{an}\subseteq V$ consisting of all these vectors has an induced
continuous $G$-action and is endowed with the subspace topology
arising from the embedding $o: V_{an}\rightarrow\Can(G,V)$. By
[E1], Prop. 2.1.26 the inclusion $C(G,K)\subseteq\Can(G,K)$ is
continuous and dualizes therefore to an algebra map
$D^c(G,K)\rightarrow\Dg$.
\begin{theo} {\bf (Schneider-Teitelbaum)}\label{peterexact}
 Let $G$ be compact and $K/\Q_p$ be finite. Suppose that $L=\Q_p$. The map
\begin{equation}\label{peterflat}
D^c(G,K)\longrightarrow\Dg\end{equation} is faithfully flat. Given
$V\in\Bang$ the representation $V_{an}$ is a strongly admissible
locally analytic representation and $V_{an}\subseteq V$ is
norm-dense. The functor $F_{\Q_p}: V\mapsto V_{an}$ between
$\Bang$ and $\Rpa$ is exact. The dual functor equals base
extension.
\end{theo}
This is [S], Thm. 4.2/3. Given the exactness statement
$V_{an}\subseteq V$ being dense is equivalent to the functor being
faithful. However, in general, the functor is not full (cf. [E2],
end of \S 3). Furthermore, if $L\neq\Q_p$, it is generally not
exact and can be zero on objects. For example (cf. [E1], \S3) let
$G=(\ol_L,+)$ and suppose $\psi:G\rightarrow K$ is $\Q_p$-linear
but not $L$-linear. The two-dimensional representation of $G$
given by the matrix $\left(%
\begin{array}{cc}
  1 & \psi \\
   & 1 \\
\end{array}%
\right)$ is an extension of the trivial representation by itself
but not locally $L$-analytic.

To study $F_L$ we have to introduce another functor. Let $G$ be an
{\it arbitrary} locally $L$-analytic group and $\Q_p\subseteq K$
be discretely valued. Given $V\in\Rpn$ we call $v\in V$ a {\it
locally $L$-analytic vector} if $o_v\in\Can(G_0,V)$ lies in the
subspace $\Can(G,V)$. Denote the space of these vectors, endowed
with the subspace topology from $V$, by $V_{an}$. Since
translation on $G$ is locally $L$-analytic $V_{an}$ has an induced
continuous $G$-action. In the following we will show that the
correspondance $V\mapsto V_{an}$ induces a functor
\[F_{\Q_p}^L: \Rpan\rightarrow\Rpa.\] Given $V\in\Rp$ the Lie
algebra $\lie$ acts on $V$ via continuous endomorphisms
\[\x v:=\frac {d}{dt}\,\exp(t\x)v|_{_{t=0}}\]
for $\x\in\lie,~v\in V$. Denote by $\Can_1(G,K)$ the local ring at
$1\in G$ as introduced before. Given $f\in\Can(G,K)$ denote its
image in $\Can_1(G,K)$ by $[f]$. Viewing $\lie_0$ as point
derivations on $\Can_1(G_0,K)$ restricting derivations to
$\Can_1(G,K)$ induces a map
$L\otimes_{\Q_p}\lie_0\rightarrow\lie$. Denote its kernel by
$\hd$. Let $(\hd)$ denote the two-sided ideal generated by $\hd$
inside $U(\lie_0,K)$ as well as in $\Dgo$. It equals the kernel of
the quotient maps $U(\lie_0,K)\rightarrow U(\lie,K)$ as well as
$\Dgo\rightarrow\Dg$ (by straighforward generalizations of [Sch],
Lem. 5.1)
\begin{lem}\label{crit}
An element $f\in\Can(G_0,K)$ is locally $L$-analytic at $1\in G$
if and only if the space of derivations $\hd$ annihilates $[f]$.
\end{lem}
\begin{pr}
The function $f$ is locally $L$-analytic at $1$ if and only if
this is true for $[f]$. By Thm. \ref{jan} $U(\lie_0,K)$ is
Fr\'echet-Stein. Hence $(\hd)\subseteq U(\lie_0,K)$ being finitely
generated is closed. By [B-TVS], IV.2.2 Corollary the natural map
$\Can_1(G_0,K)/\Can_1(G,K)\rightarrow (\hd)'$ is an isomorphism
whence the claim follows.
\end{pr}
\begin{lem}\label{h-inv}
Given $V\in\Rpn$ of compact type one has $V_{an}=V^{\hd}$ as
subspaces of $V$. In particular, $V_{an}\subseteq V$ is closed.
\end{lem}
\begin{pr}
This follows also from [E1], Prop. 3.6.19 but we give a proof in
the present language. We may assume that $G$ is compact. Suppose
first that $V=\Can(G_0,K)$. The inclusion $V_{an}\subseteq
V^{\hd}$ is clear from Lem. \ref{crit}. Let $f\in V^{\hd}$ and
$g\in G$. Denote by ${\rm Ad}$ the adjoint action of $G$. Since
$\hd$ is $L\otimes_{\Q_p}{\rm Ad}(g)$-stable the identity $g\x
g^{-1}.v={\rm Ad}(g)\x.v$ for $v\in V$ implies that $g.f$ (left
regular action) lies in $V^{\hd}$ whence is locally analytic at
$1\in G$ by Lem. \ref{crit}. This settles the case
$V=C^{an}(G_0,K)$. For general $V\in\Rpn$ of compact type
equipping $\Can(G_0,K)\hat{\otimes}_KV$ with the diagonal action
(trivial on the second factor) the topological vector space
isomorphism
$\Can(G_0,K)\hat{\otimes}_KV\stackrel{\sim}{\longrightarrow}\Can
(G_0,V)$ ([E1], Prop. 2.1.28) becomes $G$-equivariant. By
continuity one obtains $\Can(G_0,V)^{\hd}=$ closure of
$\Can(G_0,K)^{\hd}\otimes_KV$ which, by the first step, equals
$\Can(G,V)$. Hence, for $v\in V$ we have
$\lie^0\,v=0\Leftrightarrow\lie^0\, o_v=0\Leftrightarrow
o_v\in\Can(G,V)$ whence $V^{\hd}=V_{an}$.
\end{pr}
We remark that, by definition, admissible locally analytic
representations are, in particular, vector spaces of compact type.
\begin{prop}\label{comp}
The correspondance $V\mapsto V_{an}$ induces a left exact functor
\[F_{\Q_p}^L: \Rpan\rightarrow\Rpa.\]
The dual functor equals base extension.
\end{prop}
\begin{pr}
By Lem. \ref{h-inv} $V_{an}\subseteq V$ is closed and hence of
compact type. If $H\subseteq G$ is a compact open subgroup then
since $\Do/(\hd)=\Dh$ the strong dual $(V_{an})'_b$ lies in
$\C_{H_0}\cap\M (\Dh)=\C_H$. Thus, $V_{an}\in\Rpa$. It is
immediate that $V\mapsto V_{an}$ is a functor which is left exact
(Lem. \ref{h-inv}). Putting $\hdh:=\hd$ in Prop. \ref{translate}
below the last claim follows.
\end{pr}

\begin{cor}
Let $G$ be compact and $K/\Q_p$ be finite. One has
$F_L=F_{\Q_p}^L\circ F_{\Q_p}$ whence $F_L$ is left exact and the
dual functor equals base extension.
\end{cor}
\begin{pr}
Given $V\in\Bang$ the identity $F_L(V)=F_{\Q_p}^L\circ
F_{\Q_p}(V)$ as abstract $K[G]$-modules is clear from the
definitions. Since $\Can(G,V)\subseteq\Can(G_0,V)$ is a closed
topological embedding (by a straightforward generalization of
[ST4], Lem. 1.2) the topology on the space
$F_L(V)=F_{\Q_p}(V)\cap\Can(G,V)$ coincides with the one induced
by $F_{\Q_p}(V)$. By definition this equals the topology of
$F_{\Q_p}^L\circ F_{\Q_p}(V)$. The last claim follows from Thm.
\ref{peterexact} and the last proposition by associativity of the
tensor product.
\end{pr}

Remarks:

1. Since $F_L$ is not exact we obtain from the corollary that as a
rule, the map $D^c(G,K)\rightarrow\Dg$ is not flat for $L\neq\Q_p$
(but see Thm. \ref{faith}). In view of the characterization as
certain Lie invariants (Lem. \ref{h-inv}) one may ask whether
flatness holds when $G$ is semisimple. This is answered negatively
by Cor. \ref{Emerton} below.

2. On certain interesting subcategories of $\Rpan$ the functor
$F^L_{\Q_p}$ may very well be exact. To give an example recall
that $V\in\Rpa$ is called {\it locally $U(\lie)$-finite} if for
all $x\in V'_b$ the orbit $U(\lie)x$ is contained in a finite
dimensional $K$-subspace of $V'_b$. These representations are
studied in [ST1]. Let ${\rm Rep}^{a,f}_K(G)$ denote the full
abelian subcategory of $\Rpa$ consisting of these representations.
We claim: if $G$ is semisimple passage to analytic vectors is an
{\it exact} functor ${\rm Rep}^{a,f}_K(G_0)\rightarrow {\rm
Rep}^{a,f}_K(G)$. Indeed: let $V\in {\rm Rep}^{a,f}_K(G_0)$. It
suffices to see that $H^1(\lie^0,V)=0$. The $\lie^0$-module $V'_b$
is a direct limit of finite dimensional ones $W$. Since $\lie_0$
and thus $\lie^0$ are semisimple Lie algebras the first Whitehead
lemma together with Lem. \ref{isodual} for $*=1$ yields
$H_1(\lie^0,W)=0$. Since $H_1(\lie^0,.)$ commutes with direct
limits using Lem. \ref{isodual} again we obtain $H^1(\lie^0,V)=0$.

\section{$\sigma$-analytic vectors}
%%%%%%This section was changed on 04.03.2008 and in July 2008

In this section $G$ denotes a compact $d$-dimensional locally
$L$-analytic group. Under the assumption that $L/\Q_p$ is
unramified we will prove a generalization to Thm.
\ref{peterexact}.

So let us first {\bf assume} that $L/\Q_p$ is Galois. We start
with a result on the Fr\'echet-Stein structure on $U(\lie,K)$.
Given $\sigma\in Gal(L/\Q_p)$ let $G_\sigma$ be the scalar
restriction of $G$ via $\sigma: L\rightarrow L$ ([B-VAR], 5.14.1).
It is a compact locally $L$-analytic group. Denote by
$\lie_\sigma$ its Lie algebra. Of course,
$(G_\sigma)_0=G_0,~(\lie_\sigma)_0=\lie_0$ since $\sigma$ is
$\Q_p$-linear. Put $G_{L/\Q_p}:=\prod_\sigma G_\sigma$ and $\gl$
for its Lie algebra. There is a commutative diagram of locally
convex $K$-vector spaces
\[
\begin{CD}
\Can(G_{L/\Q_p},K)@>\Delta^*\circ\iota>>\Can(G_0,K)\\
@V[.]VV    @VV[.]V\\
\Can_1(G_{L/\Q_p},K)@>\Delta^*\circ\iota>>\Can_1(G_0,K)
\end{CD}
\]
where the horizontal arrows are induced functorially from the
diagonal map $\Delta: G_0\mapsto(\prod_\sigma G_\sigma)_0$
together with the canonical embedding $\iota$ (cf. (\ref{emb})).
\begin{lem}The lower horizontal map is bijective.
\end{lem}
\begin{pr}
We may assume that $G$ admits a global chart $\phi$, that $L=K$
and that ${\rm dim}_LG=1$. Using $\phi$ and the induced global
chart for $G_\sigma$ resp. $G_0$ we arrive at a map
$\Can_1(\prod_\sigma
L_\sigma,L)\rightarrow\Can_1(\Q_p^n,L)\simeq\Can_1(L^n,L)$ where
the first map depends on a choice of $\Q_p$-basis $v_i$ of $L$ and
the second identification as rings is the obvious one. Tracing
through the definitions shows that this map is induced by the
$L$-linear isomorphism $\prod_\sigma L_\sigma\simeq
L\otimes_{\Q_p}L\simeq\oplus_i Lv_i\simeq L^n$ and hence, is
bijective.
\end{pr}
Passing to strong duals yields a commutative diagramm of
topological algebras
\[
\begin{CD}
U(\lie_0,K)@>{\sim}>>
U(\gl,K)\\
@V[.]'VV    @V[.]'VV\\
D(G_0,K)@>\iota'\circ\Delta_*>>D(\Gl,K).
\end{CD}
\]
We abbreviate in the following $\varphi:=\iota'\circ\Delta_*$.

%%%%%%%%%%%%%%%%%%%%%%%%%%%%%%%%%%%%%%%%%%%%Change
\begin{lem}\label{bijective}
Suppose $G$ is uniform*. For each $r$ sufficiently close to $1$
the above diagram extends to a commutative diagram
\[
\begin{CD}
U_r(\lie_0,K)@>>>
U_r(\gl,K)\\
@VVV    @VVV\\
D_r(G_0,K)@>{\trr}>>D_r(\Gl,K)
\end{CD}
\]
of Banach algebras where the upper horizontal map is injective,
norm-decreasing and has dense image.
\end{lem}
\begin{pr}
Let $||.||_r$ be a fixed norm on $\Dgo=D((G_\sigma)_0,K)$ and let
$\nr^{\sigma}$ be the quotient norm on $\Dgs$ under
$\iota'_\sigma: D((G_\sigma)_0,K)\rightarrow\Dgs$. Let $\x_j$ and
$v_i$ be bases for the uniform* group $G$. In particular
$\x_1,...,\x_d$ is an $L$-basis for $\lie$ and $v_1=1,...,v_n$ is
a $\Z_p$-basis for $\ol_L$. The elements $h_{ij}:=\exp(v_i\x_j)$
are a minimal set of topological generators for $G_0$. Putting
$b_{ij}:=h_{ij}-1$ we obtain from (\ref{ortho}) that $\Dgo$
consists of certain series $\lambda=\sum_{\al\in\N_0^{nd}}
d_\al\db^\al$ where $||\lambda||_r=\sup_\al|d_\al|r^{\kap\al}$.
The product version of these considerations yields norms $
||.||_r^{(\sigma)}$ resp. quotient norms $\nr^{(\sigma)}$ on
$D(\prod_\sigma G_0,K)$ resp. $D(\Gl,K)$. Now $\tr$ is induced
from $\Delta$ and $D((\prod_\sigma G_\sigma)_0,K)\rightarrow
D(\prod_\sigma G_\sigma,K)$ where the second map is certainly
norm-decreasing with respect to $||.||_r^{(\sigma)}$ and
$\nr^{(\sigma)}$. Furthermore $\Delta_*(b_{ij})=\Delta(h_{ij})-1$.
Since all $nd$ elements $\Delta(h_{ij})$ are pairwise different
modulo the first step $(\prod_\sigma G_0)^p$ in the lower
$p$-series of the uniform group $\prod_\sigma G_0$ the discussion
in [DDMS], 4.2 shows that they may be completed to a minimal
ordered system of generators for $\prod_\sigma G_0$. Since the
norm $||.||_r^{(\sigma)}$ does not depend on a particular choice
of such system (cf. discussion in [ST5] after Thm. 4.10) one
obtains
$||\Delta(h_{ij})-1||_r^{(\sigma)}=||(h_{ij},1,...)-1||_r^{(\sigma)}=||h_{ij}-1||_r$
whence it easily follows that $\Delta_*$ is an isometry. Then
$\tr$ is norm-decreasing with respect to $||.||_r$ and
$\nr^{(\sigma)}$ which yields the completed diagram and its
commutativity. It is clear that the upper horizontal map is
norm-decreasing with dense image whence it remains to establish
injectivity. We first show that the inverse map
\[\tr^{-1}: U(\gl,K)\rightarrow U(\lie_0,K)\] is norm-decreasing when both sides are given suitable norm
topologies. By Thm. \ref{jan} the rings $U_r(\lie_0,K)$ resp.
$U_r(\lie_\sigma,K)$ are certain noncommutative power series rings
in the "variables" $\partial_{ij}:=v_i\x_j\in\lie_0$ resp.
$\partial_{\sigma,j}:=\x_j\in\lie_\sigma$. More precisely,
$U_r(\lie_\sigma,K)$ consists of all formal series
$\sum_{\beta\in\N_0^d}d_\beta\partial_\sigma^\beta$ where
$d_\beta\in K,~
\partial_\sigma^\beta:=\partial_{\sigma,1}^{\beta_1}\cdot\cdot\cdot\partial_{\sigma,d}^{\beta_d}$
and $||d_\beta\partial_\sigma^\beta||_{\bar{r}}^\sigma\rightarrow
0$ for $|\beta|\rightarrow\infty$. By Prop. \ref{tobe}
$\nr^\sigma$ is multiplicative, the topological $K$-basis
$\partial_\sigma^\beta$ for $U_r(\lie_\sigma,K)$
 is even orthogonal with respect to $\nr^\sigma$ and $||\partial_{\sigma,j}||_{\bar{r}}^\sigma
 =||\partial_{1j}||_r=c_r.$ Given a generic element of $U(\gl,K)$, say
$\lambda:=\sum_{\beta\in\N_0^{nd}}
d_\beta\prod_\sigma\partial_\sigma^{\beta_{\sigma.}}$ with
$d_\beta\in K$,~
$\partial_\sigma^{\beta_{\sigma.}}=\partial_{\sigma,1}^{\beta_{\sigma,1}}
\cdot\cdot\cdot\partial_{\sigma,d}^{\beta_{\sigma,d}}$ we have
$||\lambda||_{\bar{r}}^{(\sigma)}=\sup_\beta|d_\beta|(c_r)^{|\beta|}$
since the elements $\prod_\sigma\partial_\sigma^{\beta_{\sigma.}}$
are orthogonal with respect to $||.||_{\bar{r}}^{(\sigma)}$ and
this latter norm is multiplicative. After these remarks consider
the map $L\otimes_{\Q_p}L\rightarrow\prod_\sigma L_\sigma$. Let
$\sum_{a_\sigma,b_\sigma}a_\sigma\otimes b_\sigma$ be the inverse
image of $1\in L_\sigma$ where we may assume that
$b_\sigma\in\ol_L$. Choosing $b_\sigma^{(i)}\in\Z_p$ such that
$\sum_ib_\sigma^{(i)}v_i=b_\sigma$ we put
\[s:=\sup_{\sigma,a_\sigma,b_\sigma,i}|a_\sigma
b_\sigma^{(i)}|.\]
We have by definition of the map $\tr$ that
\[\tr^{-1}(\partial_{\sigma,j})=\sum_{a_\sigma,b_\sigma^{(i)},i}a_\sigma
b_\sigma^{(i)}\partial_{ij}.\] Hence,
\[||\tr^{-1}(\partial_{\sigma,j})||_r=\sup_i|\sum_{a_\sigma,b_\sigma^{(i)},i}a_\sigma
b_\sigma^{(i)}|c_r\leq sc_r\] using that the $\partial_{ij}$ are
orthogonal. Given $\lambda\in U(\gl,K)$ as above define another
norm on $U(\gl,K)$ via
\[||\lambda||_{(\bar{r})}^{(\sigma)}:=\sup_\beta|d_\beta|(sc_r)^{|\beta|}\]
and let $U_{(r)}(\gl,K)$ be the completion. We obtain
\[||\tr^{-1}(\lambda)||_r\leq\sup_\beta |d_\beta|
\prod_{\sigma,j}
||\tr^{-1}(\partial_{\sigma,j})||_r^{\beta_{\sigma,j}} =\sup_\beta
|d_\beta| (sc_r)^{|\beta|}=||\lambda||_{(\bar{r})}^{(\sigma)}\]
using that $||.||_r$ is multiplicative. Thus $\varphi^{-1}$ is
norm-decreasing with respect to the indicated norms. Now suppose
$r$ is sufficiently close to $1$. Then there exists $r'\leq r$
such that $sc_{r'}\leq c_r$ (note that $c_r\uparrow\infty$ for
$r\uparrow 1$) whence $U_r(\gl,K)\subseteq U_{(r')}(\gl,K)$. A
simple approximation argument shows that the map
\[U_r(\lie_0,K)\stackrel{\varphi_r}{\longrightarrow}U_r(\gl,K)\stackrel{\subseteq}{\longrightarrow}
U_{(r')}(\gl,K)\stackrel{\varphi^{-1}_{r'}}{\longrightarrow}U_{r'}(\lie_0,K)\]
equals the inclusion $U_r(\lie_0,K)\subseteq U_{r'}(\lie_0,K)$.
This finishes the proof.
\end{pr}
For the rest of this section we assume $L/\Q_p$ to be {\bf
unramified}.
\begin{lem}\label{unr}
Let $x\in\ol_L^\times$ be a lift of a primitive element $\bar{x}$
for the residue field extension of $L/\Q_p$. Suppose $G$ is
uniform* with a $\Z_p$-basis $v_1=1,...,v_n$ of $\ol_L$ of the
form $v_i=x^{i-1}$. In the situation of the last lemma the map
\[\varphi_r: U_r(\lie_0,K)\rightarrow U_r(\gl,K)\] is an isometry.
\end{lem}
\begin{pr}
Source and target of our map are filtered through the respective
norm and it suffices to see that the associated graded map
$\gr\varphi_r: \gr U_r(\lie_0,K)\rightarrow \gr U_r(\gl,K)$ is
injective. We use the notation of the preceding proof. Let
$X_{ij}$ resp. $X_{\sigma,j}$ be the principal symbol of
$\partial_{ij}$ resp. $\partial_{\sigma,j}$. Then $\gr
U_r(\lie_0,K)=(\gor K)[X_{11},...,X_{nd}]$ and $U_r(\gl,K)=(\gor
K)[X_{\sigma_1,1},...,X_{\sigma_n,d}]$. Denote for each $\sigma\in
Gal(L/\Q_p)$ by $\bar{\sigma}$ the induced Galois automorphism on
residue fields. Let $F_j$ be the ring homomorphism
\[F_j: (\gor K)[X_{1j},...,X_{nj}]\longrightarrow
(\gor
K)[X_{\sigma_1,j},...,X_{\sigma_n,j}],~X_{ij}\mapsto\sum_{\sigma}
\bar{\sigma}^{-1}\bar{v}_i.X_{\sigma,j}.\] Then $\gr\varphi$
equals, in the obvious sense, $F_1\otimes_{\gor K}...\otimes_{\gor
K}F_d$ whence, by induction, it suffices to prove bijectivity of
$F_1$. Now $F_1$ respects the grading by total degree whence it is
enough to prove bijectivity on homogeneous components. Since each
latter is free of finite rank over the principal ideal domain $\gor K$
%and source and target have the same rank
 it suffices to prove surjectivity in each component %since the rank is additive in short exact sequences
 %over a PID. More precisely, if F_1 is surjective in a component then rk Im=rk target
 %and therefore rk kernel=0. But if the kernel is nonzero it is free,
 %being a submodule of a free module, hence the kernel must be zero.
or, since
$F_1$ is a ring homomorphism, in the degree $1$ component.%since the algebra generators of
%the target lie in the degree one component.
 But the representing matrix $(a_{ij})_{i,j=1,...,n}$ of the
degree $1$ part of $F_1$ with respect to the $\gor K$-bases $X_i,
i=1,...,n$ resp. $X_{\sigma}, \sigma\in Gal(L/\Q_p)$ on source
resp. target has the shape
$a_{ij}=(\bar{\sigma}^{-1}_i\bar{x})^{j-1}$ and hence determinant
$\prod_{j>i}(\bar{\sigma}^{-1}_j\bar{x}-\bar{\sigma}^{-1}_i\bar{x})$.
Since $L/\Q_p$ is unramified and $\bar{x}$ generates the residue
extension this determinant is a nonzero element of the residue
field of $L$ whence a unit in $\gor K$.
\end{pr}
\begin{lem} \label{localflat}
Let $r$ be sufficiently close to $1$. The ring extension \[\trr:
\Dgor\longrightarrow D_r(\Gl,K)\] is faithfully flat.
\end{lem}
\begin{pr}
We prove only the left version (the right version follows
similarly). Suppose $H\subseteq G$ is a normal open subgroup which
is uniform*. Endow $\Dg$ with the Fr\'echet-Stein structure
induced by $\Dh$ as explained in section $2$. One obtains a
commutative diagram
\[
\begin{CD}
 \Dor@>>> D_r(\Hl,K)\\ @VVV    @VVV\\
\Dgor@>{\trr}>>D_r(\Gl,K).
\end{CD}
\]
Let $\Rep$ be a system of coset representatives for $G/H$. Choose
a system of coset representatives $\Rep'$ in $G^0$ for coset
representatives of the cokernel of the inclusion
$G/H\stackrel{\Delta}{\longrightarrow}\prod_\sigma
G_\sigma/H_\sigma=\Gl/\Hl$. Then $\trr(\Rep)\Rep'$ equals a system
of coset representatives for $\Gl/\Hl$ and so the vertical arrows
in the above diagram are finite free ring extensions on $\Rep$
resp. $\trr(\Rep)\Rep'$. Consider the map of left
$D_r(G_0,K)$-modules
\begin{equation}\label{comp2}
\oplus_{g'\in\Rep'}~\Dgor\otimes_{D_r(H_0,K)}D_r(\Hl,K)\longrightarrow
D_r(\Gl,K)\end{equation} induced by
$(\lambda\otimes\mu)_{g'\in\Rep'}\mapsto
\sum_{g'\in\Rep'}\trr(\lambda)\mu\,g'$. On the level of vector
spaces this map factores through
\[\oplus_{g\in\Rep,g'\in\Rep'}~gD_r(\Hl,K)\stackrel{\sim}{\rightarrow} \sum_{g,g'}
\trr(g)D_r(\Hl,K)g'=D_r(\Gl,K)\] and hence, is bijective (using
that $g'D_r(\Hl,K)=D_r(\Hl,K)g'$ by normality of $\Hl$ in $\Gl$).
It therefore suffices to establish the claim for $H$. In other
words, we may assume in the following that $G$ is uniform*. By the
construction of uniform* subgroups (cf. \cite{Sch}, Cor. 4.4) we
may assume that the associated $\Z_p$-basis $v_1=1,...,v_n$ is as
described in Lem. \ref{unr}. This lemma together with Lem.
\ref{bijective} then yields a commutative diagram (+)
\[
\begin{CD}
 U_r(\lie_0,K)@>>> U_r(\gl,K)\\ @VVV    @VVV\\
\Dgor@>{\trr}>>D_r(\Gl,K)
\end{CD}
\]
where the upper vertical arrow is an isometry with dense image.
Let $\dX:=\{\x_1,...,\x_d\}$ be the associated $L$-basis of
$\lie_L$ for $G$. Then the $nd$ elements $\exp(v_i\x_j)$ are a
minimal set $S$ of topological generators for $G$. Hence,
according to Thm. \ref{jan}, the left vertical arrow is finite
free on a basis $\Rep$ in $\Z[G]$. Putting
$\db^\al=b_{11}^{\al_{11}}\cdot\cdot\cdot b_{nd}^{\al_{nd}},
b_{ij}=\exp(v_i\x_j)-1\in\Z[G]$ one has
$\Rep=\{\db^\al,~\al_{ij}<\epsilon(r,p){\rm ~for~all~}ij\}$ with
$\epsilon(r,p)$ depending only on $r$ and $p$. Apply this to the
group $\Gl$ as well. More precisely, $\Gl$ is uniform* with
$v_1,...,v_n$ and $n=|Gal(L/\Q_p)|$ copies of $\dX$ as bases. By a
previous argument the set $\trr(S)\subseteq \Gl$ may be completed
to a minimal set $S'$ of generators for the uniform group
$\prod_\sigma G_0$. Choose an ordering $h'_1,...,h'_{ndn}$ of $S'$
such that $\trr(S)=\{h'_1,...,h'_{nd}\}$. Put $b'_{k}=h'_k-1$ and
form the set $\Rep':=\{\db'^\al,~\al_{k}<\epsilon(r,p){\rm
~for~all~}k\}$. Again by Thm. \ref{jan} the right vertical arrow
is finite free on the basis $\Rep'$. Let
$\Rep'_<=\{\db'^\al\in\Rep', \al_k=0 {\rm ~for~all~}k>nd\}$ and
$\Rep'_>=\{\db'^\al\in\Rep', \al_k=0 {\rm ~for~all~}k\leq nd\}$.
The chosen ordering of $S'$ implies that $\trr(\Rep)=\Rep'_<$. Now
recall that $\Dgor$ and $D_r(\Gl,K)$ are Zariski rings (cf.
section $2$) with respect to the norm filtrations and $\trr$ is
norm-preserving (Lem. \ref{bijective}) whence a filtered morphism.
Thus, it suffices to prove faithful flatness of the graded map
$\gr\trr: \gr\Dgor\rightarrow\gr D_r(\Gl,K)$ ([LVO], II.1.2.2). To
do this consider the graded version of (+) (where the upper
objects are formed with respect to the induced filtrations via the
vertical inclusions of (+))
\[
\begin{CD}
 \gr U_r(\lie_0,K)@>{\sim}>> \gr U_r(\gl,K)\\ @VVV    @VVV\\
\gr \Dgor@>{\gr\trr}>>\gr D_r(\Gl,K).
\end{CD}
\]
Recall that all rings occuring are commutative. For the
 rest of this proof $\sigma(.)$ always denotes the principal symbol map. By Prop. \ref{tobe} the
vertical arrows are finite free on the basis $\sigma(\Rep)$ resp.
$\sigma(\Rep')$ and
$\sigma(\Rep')=\{st,~s\in\sigma(\Rep'_<),~t\in\sigma(\Rep'_>)\}$.
Since $\gr\trr\neq 0$ one has $\gr\trr(\sigma(b_{ij}))\neq 0$ for
all $ij$ whence
$\gr\trr(\sigma(\Rep))=\sigma(\trr(\Rep))=\sigma(\Rep'_<)$. Then
the map of $\gr\Dgor$-modules
\[\oplus_{x'\in\sigma(\Rep'_>)}~\gr\Dgor\otimes_{\gr U_r(\lie_0,K)}\gr U_r(\gl,K)
\longrightarrow \gr D_r(\Gl,K)
\] induced by $(\lambda\otimes\mu)_{x'\in\sigma(\Rep'_>)}\mapsto \sum_{x'\in\sigma(\Rep'_>)}
\gr\trr(\lambda)\,\mu\,x'$ is bijective. It follows that $\gr
D_r(\Gl,K)$ is a finite free $\gr \Dgor$-module (of rank
$\#\Rep'_>$).
\end{pr}
\begin{cor} The map $\tr$ is injective.\end{cor}
\begin{pr}
Since all $\trr$ are injective (for $r$ sufficiently close to $1$)
and compatible with transition with respect to $r'\leq r$ the map
$\tr$ is injective by left-exactness of the projective limit.
\end{pr}
%%%%End of change
For the rest of this section we assume $K/\Q_p$ to be finite.
Consider the faithfully flat algebra map $D^c(G,K)\rightarrow\Dgo$
stated in (\ref{peterflat}). We compose it with $\tr$ and show our
first main result.
\begin{theo}\label{faith}
The map $D^c(G,K)\longrightarrow D(\Gl,K)$ is faithfully flat.
\end{theo}
\begin{pr}
We show only left faithful flatness. For flatness we are reduced,
by the usual argument, to show that the map
$D(\Gl,K)\otimes_{D^c(G,K)}J\longrightarrow D(\Gl,K)$ is injective
for any left ideal $J\subseteq D^c(G,K)$. The ring $D^c(G,K)$
being noetherian the left hand side is a coadmissible
$D(\Gl,K)$-module. By left exactness of the projective limit we
are thus reduced to show that
$D_r(\Gl,K)\otimes_{D^c(G,K)}J\rightarrow D_r(\Gl,K)$ is injective
for all $r$ (sufficiently close to $1$). This is clear since
\[D^c(G,K)\rightarrow D(G_0,K)\rightarrow
D_r(G_0,K)\stackrel{\tr_r}{\longrightarrow}D_r(\Gl,K)\] is flat
(the second map by [ST5], remark 3.2).

For faithful flatness we have to show
$D(\Gl,K)\otimes_{D^c(G,K)}M\neq 0$ for any nonzero left
$D^c(G,K)$-module $M$. By the first step we may assume that $M$ is
finitely generated. Then $D(\Gl,K)\otimes_{D^c(G,K)}M$ is
coadmissible whence we are reduced, by the equivalence of
categories between coadmissible modules and coherent sheafs
([ST5], Cor. 3.3), to find an index $r$ such that
\[D_r(\Gl,K)\otimes_{D^c(G,K)}M\neq 0.\] Put
$N:=\Dgo\otimes_{D^c(G,K)}M\in\C_{G_0}$. Then $N\neq 0$ whence
$N_r\neq 0$ for some $r$ (sufficiently close to $1$). It follows
that $D_r(\Gl,K)\otimes_{D^c(G,K)}M$ equals
\[D_r(\Gl,K)\otimes_{\Dgo}N=D_r(\Gl,K)\otimes_{\Dgor}N_r\] and the
right-hand side is nonzero by faithful flatness of $\tr_r$.
\end{pr}
Each choice $\sigma\in Gal(L/\Q_p)$ gives rise to the locally
$L$-analytic manifold $K_{\sigma}$ arising from restriction of
scalars via $\sigma$. The space $\Can
(G,K_\sigma)=\Can(G_{\sigma^{-1}},K)$ is called the space of {\it
locally} $\sigma$-{\it analytic} functions ([B-VAR], 5.14.3). This
motivates the following definition: given $V\in\Rp$ an element
$v\in V$ is called a {\it locally $\sigma$-analytic vector} if
$o_v\in \Can(G_{\sigma^{-1}},V)$. Let $V_{\sigma-an}$ denote the
subspace of all these vectors in $V$.

Consider the Lie algebra map
$L\otimes_{\Q_p}\lie_0\simeq\prod_\sigma\lie_\sigma\rightarrow\lie_{\sigma^{-1}}.$
The kernel $\hs$ acts on $V$ whence one deduces as in case
$\sigma=id$ that $V_{\sigma-an}=V^{\hs}$, functorial in $V$ and
that passage to $\sigma$-analytic vectors is a left exact functor
$\Rpan\rightarrow {\rm Rep}^{a}_K(G_{\sigma^{-1}})$. Note also
that given $\sigma\neq\tau$ one has $V_{\sigma-an}\cap
V_{\tau-an}=V^{\lie_0}=V^\infty$, the smooth vectors in $V$.

We consider the base extension functor $M\mapsto
M\otimes_{D^c(G,K)}D(\Gl,K)$ on finitely generated
$D^c(G,K)$-modules and pull back to representations. This yields a
functor
\[F: \Bang\longrightarrow{\rm Rep}^{a}_K(\Gl)\]
which is exact and faithful according to Thm. \ref{faith}. Given
$V\in\Bang$ note that, by exactness and since $D^c(G,K)$ is
noetherian, the coadmissible module associated to $F(V)$ is even
finitely presented. We deduce the second main result.
\begin{theo}\label{new}
The functor $F$ is exact and faithful. Given $V\in\Bang$ the
representation $F(V)$ is strongly admissible. Viewed as a
$G_{\sigma^{-1}}$-representation, $\sigma\in Gal(L/\Q_p)$, it
contains $V_{\sigma-an}$ as a closed subrepresentation and
functorial in $V$. In case $L=\Q_p$ the functor coincides with
$F_{\Q_p}$.
\end{theo}
\begin{pr}
It remains to see the latter statements. The projection
$pr_\sigma: \prod_\sigma G_\sigma\rightarrow G_\sigma$ induces a
continuous inclusion $pr_\sigma^*:
\Can(G_\sigma,K)\rightarrow\Can(\Gl,K)$. By definition of the
locally convex topologies on both sides it is a compact locally
convex inductive limit of isometries and so, according to [E1],
Prop. 1.1.41, a topological embedding with closed image. Dualizing
we obtain a continuous algebra surjection $(pr_\sigma)_*:
D(\Gl,K)\rightarrow D(G_\sigma,K)$ exhibiting $D(G_\sigma,K)$ as
coadmissible $D(\Gl,K)$-module. It follows from [ST5], Lem. 3.8
and its proof that $V'_b\otimes_{D^c(G,K)}
D(G_\sigma)\in\C_{G_\sigma}$ lies in $\C_{\Gl}$ and that the two
canonical topologies coincide. The $D(\Gl,K)$-linear surjection
\[{\rm id}\otimes (pr_\sigma)_*: V'_b \otimes_{D^c(G,K)} D(\Gl,K)\longrightarrow V'_b\otimes_{D^c(G,K)}
D(G_\sigma,K)\] then lies in $\C_{\Gl}$ and is therefore
continuous and strict. It is also $D(G_\sigma,K)$-linear and the
right-hand side equals $(V_{\sigma^{-1}-an})_b'$ according to the
$\sigma^{-1}$-analytic version of Prop. \ref{comp} (cf. remarks
above). Passing to strong duals yields a closed
$G_{\sigma}$-equivariant topological embedding
$V_{\sigma^{-1}-an}\rightarrow F(V)$. It is natural in $V$ and
clearly onto in case $L=\Q_p$.
\end{pr}
Remark: Let $H\subseteq G$ be a compact subgroup and denote the
versions of the functor $F$ relative to $H$ resp. $G$ by $F_H$
resp. $F_G$. The natural map
\[D^c(G,K)\otimes_{D^c(H,K)}D(\Hl,K)\rightarrow D(\Gl,K)\] being an
isomorphism of bimodules is equivalent to $L=\Q_p$ (cf.
(\ref{comp2})). It follows that in case $L\neq\Q_p$ the functors
$F_H$ and $F_G$ do not commute with the restriction functors
induced by $H\subseteq G$ resp. $\Hl\subseteq \Gl$. Therefore
there is no naive generalization of the functor $F$ to noncompact
groups in the case $L\neq\Q_p$ and we leave this matter as an open
problem.

\section{Standard resolutions}
Turning back to a general extension $L/\Q_p$ (not necessarily
Galois) the functors $F^L_{\Q_p}$ and $F_L$ defined previously
remain interesting in themselves. We begin their study with some
general analysis of certain base extension functors between
coadmissible modules.

Let $G$ be a locally $L$-analytic group. For the rest of this
section we {\it fix} an ideal $\hdh$ of $L\otimes_{\Q_p}\lie_0$
stable under $L\otimes {\rm Ad}(g)$ for all $g\in G$ where Ad
refers to the adjoint action of $G$. Denote by $(\hdh)$ the
two-sided ideal generated by $\hdh$ in $L\otimes_{\Q_p} U(\lie_0)$
as well as in $\Dgo$. Put
$D:=\Dgo/(\hdh),~\C_D:=\C_{G_0}\cap\mathcal{M}(D).$ Then $\C_D$
(with $D$-linear maps) is an abelian category. If $G$ is compact
then since $(\hdh)$ is closed $D$ is a Fr\'echet-Stein algebra and
if $(\hdh)_r$ denotes the closure $(\hdh)\subseteq\Dgor$, the
coherent sheaf associated to $D$ equals $D_r:=\Dgor/(\hdh)_r$
([ST5], Prop. 3.7). If we base extend the standard complex
$U(\hdh)\otimes_L\dot{\bigwedge}\hdh$ via $U(\hdh)\subseteq\Dgo$
then the complex
\begin{equation}\label{standard}
\Dgo\otimes_L\bigwedge^l\hdh\rightarrow...\rightarrow\Dgo\otimes_L\bigwedge^0\hdh\rightarrow\Dd,\end{equation}
$l:={\rm dim}_L\hdh$, consists of finite free left $\Dgo$-modules.
\begin{lem}\label{left}\label{rs} The complex (\ref{standard})
is a free resolution of the left $\Dgo$-module $D$ by
$\Dgo$-bimodules.
\end{lem}
\begin{pr}
Let us prove that $\Dgo\otimes_L\dot{\bigwedge}\hdh$ is acyclic.
Choosing $H\subseteq G$ compact open and using $\hdh$-invariant
decompositions
\[\Dgo=\oplus_{g\in G/H}g\Do,~~~ \Dd=\oplus_{g\in G/H}g
(D(H_0,K)/(\hdh))\] we are reduced to $G$ compact. The complex
(\ref{standard}) consists then of coadmissible left $\Dgo$-modules
whence acyclicity may be tested on coherent sheafs. It thus
suffices to see that $\Dgor\otimes_L\dot{\bigwedge}\hdh$ is exact
for a fixed radius $r$. The maps $U(\lie_0,K)\rightarrow
U_r(\lie_0,K)\rightarrow\Dgor$ are flat, the first by [ST5],
remark 3.2 and the second by Thm. \ref{jan}. We are thus reduced
to show that base extending the standard complex via
$U(\hdh)\subseteq U(\hdh,K)\subseteq
U(L\otimes_{\Q_p}\lie_0,K)=U(\lie_0,K)$ is an exact operation. By
Lem. \ref{acyclic} this holds for the extension $U(\hdh)\subseteq
U(\hdh,K)$ and $U(\hdh,K)\subseteq
U(L\otimes_{\Q_p}\lie_0,K)=U(\lie_0,K)$ is clearly a free ring
extension. This proves acyclicity. Using that $\hdh$ is Ad-stable
we may endow the complex (\ref{standard}) with a right
$\Dgo$-module structure as follows. The adjoint action of $G_0$ on
$\hdh$ is locally analytic and extends functorially to a
continuous right action on $\bigwedge^q\hdh$ given explicitly via
$(\x_1\wedge...\wedge\x_q)g={\rm Ad}(g^{-1})\x_1\wedge...\wedge
{\rm Ad}(g^{-1})\x_q.$ Letting $G_0$ act on $\Dgo$ by right
multiplication we give $\Dgo\otimes\bigwedge^q\hdh$ the right
diagonal $G_0$-action which extends to a separately continuous
right $\Dgo$-module structure (cf. [ST6], Appendix). The identity
$g\x g^{-1}={\rm Ad}(g)\x$ in $\Dgo$ implies that the differential
$\partial$ of (\ref{standard}) respects the diagonal right
$G_0$-action. Since $K[G_0]\subseteq\Dgo$ is dense ([ST2], Lem.
3.1) $\partial$ respects the right $\Dgo$-module structure by
continuity.
\end{pr}

\begin{prop}\label{preservemod}\label{comparelie}
Given $X\in\Cgo$ one has ${\rm Tor}^{\Dgo}_*(X,\Dd)\in\C_D$.
Furthermore, ${\rm Tor}^{\Dgo}_*(X,\Dd)\simeq H_*(\hdh,X)$ in
$Vec_K$ natural in $X$ which is topological with respect to the
canonical topology on the left hand side. In particular, ${\rm
Tor}^{\Dgo}_*(.,\Dd)=0$ in degrees $*>{\rm dim}_L\hdh$.
\end{prop}
\begin{pr}
Let $X\in\Cgo$ be given. By the above lemma ${\rm
Tor}^{\Dgo}_*(X,\Dd)\simeq h_*(X\otimes_L\dot{\bigwedge}\hdh)$ in
$Vec_K$. By the usual argument with double complexes, the right
module structure on $X\otimes_L\dot{\bigwedge}\hdh$ makes the
isomorphism right $\Dgo$-equivariant. Now $X$ is coadmissible and
$\bigwedge^q\hdh$ is finite dimensional. Hence dualising [E1],
Prop. 6.1.5 yields that each right module
$X\otimes_L\bigwedge^q\hdh$ is coadmissible and its canonical
topology coincides with the tensor product topology. Since $\Cgo$
is abelian we obtain ${\rm
Tor}^{\Dgo}_*(X,\Dd)\in\C_{G_0}\cap\mathcal{M}(\Dd)=\mathcal{C}_D$.
The remaining statements are now clear.
\end{pr}
In the compact case the associated coherent sheafs are easily
computed.
\begin{cor}\label{sheaf}
Let $G$ be compact. Given $X\in\Cgo$ then ${\rm
Tor}^{\Dgo}_*(X,\Dd)_r={\rm Tor}^{\Dgor}_*(X_r,\Dd_r)$. Also,
${\rm Tor}^{\Dgor}_*(X_r,\Ddr)\simeq H_*(\hdh,X_r)$ in $Vec_K$.
\end{cor}
\begin{pr}
Let $P.\rightarrow X$ be a projective resolution in $\Mgn$. By
flatness of $\Dgo\rightarrow\Dgor$ the complex
$P.\otimes_{\Dgo}\Dgor\rightarrow X_r$ is a projective resolution
of the $\Dgor$-module $X_r$. Since $\Dd_r\simeq
\Dd\otimes_{\Dgo}\Dgor)$ as $(\Dgo,\Dgor)$-bimodules ([ST5], Cor.
3.1) we have as right $\Dgor$-modules
\[\begin{array}{ccl}
  {\rm Tor}^{\Dgor}_*(X_r,\Ddr) & \simeq & h_*(P.\otimes_{\Dgo}
  \Ddr) \\
  & &  \\
   & \simeq & h_*(P.\otimes_{\Dgo}\Dd)\otimes_{\Dgo}\Dgor\\& &  \\
    & \simeq & {\rm
Tor}^{\Dgo}_*(X,\Dd)\otimes_{\Dgo}\Dgor.
\end{array}\]
\end{pr}
Since $\C_{G_0}\subseteq\Mgn$ is a full embedding (and similarly
for $\C_D)$ we have
\begin{cor}\label{leftsatellite}
The functors ${\rm Tor}^{\Dgo}_*(.,\Dd)$ form a homological
$\delta$-functor between $\C_{G_0}$ and $\C_D$.
\end{cor}
\begin{prop}\label{translate} There is a commutative (up to natural isomorphism) diagram of
functors
\[
\begin{CD}
\Rpan @>{V\mapsto V^\hdh}>>\Rpan\\
@V{V\mapsto V'_b}VV    @VV{V\mapsto V'_b}V\\
\C_{G_0}@>{M\mapsto M\otimes_{\Dgo}D}>>\C_{G_0}
\end{CD}
\]
\end{prop}
\begin{pr}
Suppose first that $G$ is compact. By continuity of the Lie action
and since $\hdh$ is Ad-stable $V\mapsto V^\hdh$ is an auto-functor
of $\Rpan$. The complex $V'_b\otimes_L\dot{\bigwedge}\hdh$
consists of coadmissible right modules whence the differential is
strict. Hence  Lem. \ref{isodual} for $*=0$ implies that
restriction of functionals yields a $G_0$-isomorphism $
V_b'/V_b'\hdh\simeq (V^{\hdh})'_b$ of topological vector spaces,
functorial in $V$. By local analyticity this extends to an
isomorphism of right $\Dgo$-modules $V'_b\otimes_{\Dgo}\Dd\simeq
(V^{\hdh})'_b$ natural in $V$. This settles the compact case. If
$G$ is arbitrary the result follows easily from the compact case
by choosing a compact open subgroup.
\end{pr}

\section{Higher analytic vectors}
We are mainly interested in the choice $\hdh:=\hd$ where, as
before, $\lie^0=\ker\, (L\otimes_{\Q_p}\lie_0\rightarrow\lie)$.
Hence $D=\Dgo/(\hd)\simeq\Dg$ and $\C_D=\C_{G_0}\cap\Mg=\C_G$.
\begin{theo}\label{main}
Passage to analytic vectors $F_{\Q_p}^L$ extends to a
cohomological $\delta$-functor $(\Ti)_{i\geq 0}$ with $\Ti=0$ for
$i>([L:\Q_p]-1)\,{\rm dim}_LG$.
\end{theo}
\begin{pr}
We may clearly replace the right upper resp. lower corner in Prop.
\ref{translate} by $\Rpa$ resp. $\C_G$ without changing the
statement. Both vertical arrows are anti-equivalences between
abelian categories and therefore exact functors. By direct
calculation pulling back the functors ${\rm
Tor}^{\Dgo}_*(.,D(G,K)): \C_{G_0}\rightarrow\C_G$ (cf. Cor.
\ref{leftsatellite}) yields a cohomological $\delta$-functor
extending $F_{\Q_p}^L$. Finally, ${\rm
dim}_L\lie^0=([L:\Q_p]-1)\,{\rm dim}_LG$.
\end{pr}
The functors $\Ti$ can be expressed without referring to
coadmissible modules. Endowing $V\in\Rpan$ with the uniquely
determined separately continuous left $\Dgo$-module structure we
may consider the $G$-representation (*)
\[{\rm Ext}^i_{\Dgo}(\Dg,V)\] where the left $G$-action
comes from right multiplication on $\Dg$.
\begin{cor} The $G$-representation $(*)$ lies in $\Rpa$ and there is a natural isomorphism
\[\Ti(V)\simeq {\rm Ext}^i_{\Dgo}(\Dg,V)\] of admissible
$G$-representations.
\end{cor}
\begin{pr}
Taking cohomology on the complex ${\rm
Hom}_{\Dgo}(\Dgo\otimes_L\dot{\bigwedge}\hd,V)={\rm
Hom}_{L}(\dot{\bigwedge}\hd,V)$ yields an isomorphism ${\rm
Ext}^*_{\Dgo}(\Dg,V)\simeq H^{*}(\hd,V)$ in $Vec_K$ natural in
$V$. We endow the Ext group with the locally convex topology of
compact type of the right-hand side. Using Lem. \ref{isodual} as
well as Prop. \ref{comparelie} we obtain a natural isomorphism in
$Vec_K$
\[{\rm Tor}^{\Dgo}_*(V_b',\Dg)'_b\simeq {\rm
Ext}^*_{\Dgo}(\Dg,V)\] which is topological. To check that it is
$G$-equivariant amounts to check the $G$-equivariance of the
isomorphisms of complexes
$(h_*(V'\otimes_L\dot{\bigwedge}\hd))'\simeq
h^*((V'\otimes_L\dot{\bigwedge}\hd)')$ appearing in the proof of
Lem. \ref{isodual} and $(V'\otimes_L\dot{\bigwedge}\hd)'\simeq
{\rm Hom}_L(\dot{\bigwedge}\hd,V)$ (remark before [loc.cit.]).
These are direct computations.\end{pr}

Let $G$ be compact and $K/\Q_p$ be finite. Using the notation of
section 3 recall that $F_{\Q_p}$ is exact and preserves injective
objects (Thm. \ref{peterexact}). Hence we have the following
direct application of our results to Banach space representations.

\begin{prop}
Let $G$ be compact and $K/\Q_p$ be finite. Then
\[R^{i}F_L=\Ti\circ F_{\Q_p}{\rm~and~}R^{i}F_L=0,~i>([L:\Q_p]-1)\,{\rm dim}_LG\] for the right-derived functors of
$F_L$.
\end{prop}

We conclude with two further applications when varying $\hdh$.

1. $\sigma$-analytic representations. Assume $L/\Q_p$ is Galois.
Letting $\hdh:=\hs$ in Cor. \ref{leftsatellite} and pulling back
to representations yields: for each $\sigma\in Gal(L/\Q_p)$ the
left-exact functor $\Rpan\rightarrow{\rm
Rep}_K^{a}(G_{\sigma^{-1}}),~V\mapsto V_{\sigma-an}$ extends to a
$\delta$-functor. The higher functors vanish in degrees
$>([L:\Q_p]-1)\,{\rm dim}_LG$.

2. Smooth representations. Let ${\rm
Rep}_K^{\infty,a}(G)\subseteq\Rpan$ denote the full abelian
subcategory of smooth-admissible representations ([ST5], \S6). The
equivalence $\Rpan\simeq\C_{G_0}$ induces an equivalence ${\rm
Rep}_K^{\infty,a}(G)\simeq\C_\infty$ where
$\C_\infty=\C_{G_0}\cap\mathcal{M}(D^\infty(G,K))$ and
$D^\infty(G,K)=\Dgo/(\lie_0)$ denotes the algebra of smooth
distributions ([ST6], \S1). Since $V^\infty=H^0(\lie_0,V)$ we may
put $\hdh:=L\otimes_{\Q_p}\lie_0$, apply the results of the
preceding section and obtain: passage to smooth vectors
$\Rpan\rightarrow {\rm Rep}_K^{\infty,a} (G),~V\mapsto V^\infty$
extends to a $\delta$-functor vanishing in degrees
$>[L:\Q_p]\,{\rm dim}_LG$.

\section{Analytic vectors in induced representations}

As an application we study the interaction of the functors $\Ti$
with locally analytic induction. This implies an explicit formula
for the higher analytic vectors in principal series
representations. As usual $G$ denotes a locally $L$-analytic
group.

We let $P\subseteq G$ be a closed subgroup with Lie algebra $\p$.
Let ${\rm Ind}_P^G$ denote the locally analytic induction viewed
as a functor from admissible $P$-representations $(W,\rho)$,
finite dimensional over $K$, to admissible $G$-representations.
Explicitly,
\[{\rm
Ind}_P^G(W):=\{f\in\Can(G,W),~f(gb)=\rho(b)^{-1}f(g){\rm~for~all~}g\in
G,b\in P\}\] and $G$ acts by left translations. One has the
isomorphism of right $\Dg$-modules
\begin{equation}\label{nocheiniso}
W'_b\otimes_{D(P,K)}\Dg\stackrel{\sim}{\longrightarrow}({\rm
Ind}_{P}^{G}W)'_b\end{equation} mapping $\phi\otimes\lambda$ to
the functional $f\mapsto\lambda(g\mapsto\phi(f(g)))$ ([OS], 2.4).

Recall that $\lie^0:={\rm ker}
~(L\otimes_{\Q_p}\lie_0\rightarrow\lie)$. We {\bf assume} for the
rest of this section that there is a compact open subgroup
$G'\subseteq G$ such that an "Iwasawa decomposition"
\begin{equation}\label{iwasawa} G=G'P\end{equation} holds.
\begin{lem}\label{6.1I}
Given $X\in\Mgn$ there is a natural isomorphism
\[{\rm Tor}_*^{\Dgos}(X,D(G',K))\simeq {\rm
Tor}_*^{\Dgo}(X,D(G,K))\] as right $D(G',K)$-modules.
\end{lem}
\begin{pr}
Let $*=0$. For $X=\Dgo$ the claim follows from the fact that
$G'\subseteq G$ is open whence $\hd\subseteq\Dgos$ with
$\Dgss=\Dgos/(\hd)$. The case of arbitrary $X$ follows from this.
In general a projective resolution $P.\rightarrow X$ of $X$ as
$\Dgo$-module remains a projective resolution of $X$ as
$\Dgos$-module since $\Dgo$ is free over $\Dgos$. The claim then
follows from the case $*=0$.
\end{pr}
Now put $P':=P\cap G'$. Using (\ref{iwasawa}) restriction of
functions induces a topological $G'$-isomorphism
\begin{equation}\label{Fea}
{\rm Ind}_{P}^{G}W\stackrel{\sim}{\longrightarrow} {\rm
Ind}_{P'}^{G'}W\end{equation} where the right-hand side has the
obvious meaning ([Fea], 4.1.4).
\begin{lem}\label{6.1II}\label{6.1III} {\rm (i)} For all $Y\in \mathcal{M}(D(P,K))$ the natural map
\[Y\otimes_{D(P',K)}\Dgss\rightarrow Y\otimes_{D(P,K)}\Dg\] is
an isomorphism of $\Dgss$-modules.

{\rm (ii)} We have a commutative diagram
\[
\begin{CD}
W'\otimes_{D(P',K)}D(G',K)@>>>({\rm
Ind}_{P'}^{G'}W)'\\ @VVV    @VVV\\
W'\otimes_{D(P,K)}D(G,K)@>>>({\rm Ind}_{P}^{G}W)'
\end{CD}
\]
of right $\Dgss$-modules in which all four maps are isomorphisms.
\end{lem}
\begin{pr}
This is a straightforward generalization of [ST6], Lem. 6.1 using
(\ref{iwasawa}).
\end{pr}
Recall from last section that $F_{\Q_p}^L$ extends to a
$\delta$-functor $(R^{i}F_{\Q_p}^L)_{i\geq 0}$ (Thm. \ref{main}).
Assume that we are given a finite dimensional locally
$\Q_p$-analytic $P$-representation $W$. We abbreviate in the
following \[V:={\rm Ind}_{P_0}^{G_0}W\in\Rpan\] and study the
admissible $G$-representations $R^{i}F_{\Q_p}^L(V)$. Let
$Q.=D(P_0,K)\otimes_L\dot{\bigwedge}\p^0$ resp.
$\tilde{Q}.=\Dgo\otimes_L\dot{\bigwedge}\lie^0$ denote the
standard resolutions for the bimodules $D(P,K)$ resp. $D(G,K)$ as
referred to in Lem. \ref{rs}. We have the natural morphism of
complexes of $D(P_0,K)$-bimodules $Q.\rightarrow \tilde{Q}.$
induced by $D(P_0,K)\rightarrow\Dgo$ and $\p^0\rightarrow\lie^0$.
Tensoring with the right $D(P_0,K)$-equivariant map $W'\rightarrow
V', w\mapsto w\otimes 1$ arising from (\ref{nocheiniso}) gives
rise to a morphism of complexes of right $D(P_0,K)$-modules
\[W'\otimes_{D(P_0,K)}Q.
\rightarrow V'\otimes_{\Dgo}\tilde{Q}.\] Taking homology and
extending scalars yields a map
\[
f: {\rm Tor}_*^{D(P_0)}(W',D(P))\otimes_{D(P)}D(G)\longrightarrow
{\rm Tor}_*^{D(G_0)}(V',D(G))\] of right $D(G)$-modules where we
have abbreviated $D(G):=\Dg, D(P):=\Dp$ etc.
\begin{prop}\label{isospeed}
The map $f$ is an isomorphism.
\end{prop}
\begin{pr}
First note that the right-hand side is a priori coadmissible by
Prop. \ref{preservemod}. We now have bijective maps of right
$D(G')$-modules
\[{\rm
Tor}_*^{D(P_0)}(W',D(P))\otimes_{D(P)}D(G)\stackrel{\sim}{\rightarrow}
{\rm Tor}_*^{D(P_0)}(W',D(P))\otimes_{D(P')}D(G')
\]
(Lem. \ref{6.1III}~(i) applied to $Y:={\rm
Tor}_*^{D(P_0)}(W',D(P))$) and \[{\rm
Tor}_*^{D(P_0)}(W',D(P))\otimes_{D(P')}D(G')
\stackrel{\sim}{\rightarrow} {\rm
Tor}_*^{D(P'_0)}(W',D(P'))\otimes_{D(P')}D(G')\] (Lem. \ref{6.1I}
applied to $P'\subseteq P$ and $W'\in \mathcal{M}(\Dpo)$. Their
composite fits into the diagram of right $D(G')$-modules
\[
\begin{CD}
{\rm Tor}_*^{D(P_0)}(W',D(P))\otimes_{D(P)}D(G) @>f>>{\rm
Tor}_*^{D(G_0)}(V',D(G))\\ @VVV    @VVV\\{\rm
Tor}_*^{D(P_0')}(W',D(P'))\otimes_{D(P')}D(G') @>>> {\rm
Tor}_*^{D(G_0')}(V',D(G'))
\end{CD}
\]
where the right hand vertical arrow is due to Lem. \ref{6.1I} and
bijective. The lower horizontal arrow is defined analogously to
the upper one using Lem. \ref{6.1II}~(ii). Tracing through the
definitions of the maps involved this diagram commutes. We may
thus assume that $G$ is compact. Then both sides of our map are
coadmissible: ${\rm Tor}^{D(P_0)}_*(W',D(P))$ is finite
dimensional over $K$ (Prop. \ref{preservemod}) hence is a finitely
presented $D(P)$-module. We introduce another map of right
$D(G)$-modules
\begin{equation}\label{help2}f': {\rm
Tor}_*^{D(P_0)}(W',D(P))\otimes_{D(P)}D(G)\longrightarrow {\rm
Tor}_*^{D(G_0)}(V',D(G))\end{equation} as follows. Choose a
projective resolution $P.\rightarrow W'$ by right $D(P_0)$-modules
according to Lem. \ref{help} below. Then
$\tilde{P}.\,:=P.\otimes_{D(P_0)}D(G_0)$ is a projective
resolution for $W'\otimes_{D(P_0)} D(G_0)=V'$ whence the natural
map
\[P.\otimes_{D(P_0)} D(P)\rightarrow
\tilde{P}.\otimes_{D(G_0)}D(G)\] induced by
$m\otimes\lambda\mapsto m\otimes 1\otimes\lambda$ for $m\in
P_n,~m\otimes 1\in\tilde{P}_n,~\lambda\in D(P)\subseteq D(G)$
gives our map $f'$. We claim that it is bijective: by
coadmissibility this may be tested on coherent sheafs. Let us
realize $D(G_0)$ and $D(G)$ as Fr\'echet-Stein algebras via the
families of norms appearing in Prop. \ref{family}. In particular
$D_r(P)\rightarrow D_r(G)$ is flat for all $r$. Denote by $W'_r$
resp. $V'_r=W'_r\otimes_{D_r(P_0)}D_r(G_0)$ the coherent sheafs
associated to $W'$ resp. $V'$. Then the coherent sheafs associated
to both sides of (\ref{help2}) are given by ${\rm
Tor}_*^{D_r(P_0)}(W_r',D_r(P))\otimes_{D_r(P)}D_r(G)$ resp. ${\rm
Tor}_*^{D_r(G_0)}(V_r',D_r(G))$ according to Cor. \ref{sheaf}. Put
$P_r.:=P.\otimes_{D(P_0)}D_r(P_0)$ resp.
$\tilde{P}_r.:=\tilde{P}.\otimes_{D(G_0)}D_r(G_0)$. By [ST5],
remark 3.2 $P_r.\rightarrow W'_r$ resp. $\tilde{P}_r.\rightarrow
V'_r$ are projective resolutions of $W'_r$ resp. $V'_r$ and the
map
\[
f'\otimes_{D(G)}D_r(G): {\rm
Tor}_*^{D_r(P_0)}(W_r',D_r(P))\otimes_{D_r(P)}D_r(G)\rightarrow
{\rm Tor}_*^{D_r(G_0)}(V_r',D_r(G))\] coincides with the one
induced by
\[P_r.\otimes_{D_r(P_0)} D_r(P)\rightarrow
\tilde{P}_r.\otimes_{D_r(G_0)}D_r(G).\] Since $D_r(P)\rightarrow
D_r(G)$ is flat $f'\otimes_{D(G)} D_r(G)$ is bijective and since
this holds for all $r$ the map $f'$ is bijective. Using a standard
double complex argument now shows that $f=f'$ whence the
proposition.
\end{pr}
The following lemma was used in the preceding proof.
\begin{lem}\label{help}
Assume $G$ is compact. There is a projective resolution
$P.\rightarrow W'$ by right $D(P_0)$-modules such
$P.\otimes_{D(P_0)}D(G_0)$ is acyclic.
\end{lem}
\begin{pr}
By density of analytic vectors (Thm. \ref{peterexact}) $W'$ being
finite dimensional over $K$ implies that $W'\otimes_{D^c(P_0)}
D(P_0)=W'$. Now choose a finite free resolution $P'.$ of the
$D^c(P_0)$-module $W'$. By the flatness result (\ref{peterflat})
$P.\,:=P'.\otimes_{D^c(P_0)} D(P_0)$ is a finite free resolution
of the $D(P_0)$-module $W'$ and it remains to see the last
statement. Now every kernel $K_n\subseteq P_n$ of the differential
in $P.$ is a finitely presented $D(P_0)$-module whence the
morphism $K_n\otimes_{D(P_0)} D(G_0)\rightarrow
P_n\otimes_{D(P_0)} D(G_0)$ lies in $C_{G_0}$. Its injectivity
follows therefore on coherent sheafs from flatness of
$D_r(P_0)\rightarrow D_r(G_0)$ and left-exactness of the
projective limit using the Fr\'echet-Stein structure of Prop.
\ref{family}.
\end{pr}
\begin{theo}\label{interaction}
The functors $R^{i}F_{\Q_p}^L$ commute with induction: given a
finite dimensional locally $\Q_p$-analytic $P$-representation $W$
one has an isomorphism
\[R^{i}F_{\Q_p}^L\circ{\rm Ind}_{P_0}^{G_0}(W)\simeq {\rm
Ind}_{P}^{G}\circ \Ti(W)\] as admissible $G$-representations
functorial in $W$.
\end{theo}
\begin{pr}
This follows from dualising the isomorphism in the preceding
proposition which is, by construction, functorial in $W$.
\end{pr}
The functor ${\rm Ind}_P^G$ is nonzero on objects ([Fea], Satz.
4.3.1) whence
\begin{cor}\label{vertausch}
We have $\Ti({\rm Ind}_{P_0}^{G_0}W)\neq 0$ if and only if
$\Ti(W)\neq 0$. In particular $\Ti({\rm Ind}_{P_0}^{G_0}W)=0$ for
all $i>([L:\Q_p]-1)\,{\rm dim}_L\mathfrak{p}$.
\end{cor}
The above results apply in particular when $G$ equals the
$L$-points of a connected reductive group over $L$ and $P\subseteq
G$ is a parabolic subgroup. If $W$ is a one dimensional
$P$-representation $K_\chi$ given by a locally $\Q_p$-analytic
character $\chi: P\rightarrow K^\times$ we may determine the
vector space $\Ti(K_\chi)$ completely. For simplicity we assume
that $G$ is quasi-split and let $P:=B$ be a Borel subgroup with
Lie algebra $\bo$. Then $\bo=\tor\un$ (semidirect product) where
$\tor$ is a maximal toral subalgebra and $\un=[\bo,\bo]$. Let
$\bo^0$ be the kernel of $K\otimes_{\Q_p}\bo_0\rightarrow
K\otimes_L\bo$ and define $\tor^0$ and $\un^0$ analogously. Then
$K_\chi$ is a $\bo^0$-module via $K\otimes_{\Q_p}{\rm d}\chi$
where ${\rm d}\chi: \bo_0\rightarrow K$ denotes the differential
of $\chi$. We denote this module as well as the induced (note that
$[\bo^0,\bo^0]=\un^0)$ $\tor^0$-module by $K_{{\rm d}\chi}$.
\begin{cor}
There is an isomorphism in $Vec_K$
\[\Ti(K_\chi)\simeq\sum_{j+k=i}\,\bigwedge^{j}\tor^0\otimes_K
H^k(\un^0,K_{{\rm d}\chi})^{\tor^0}.\]
\end{cor}
\begin{pr}
By Prop. \ref{comparelie} we have that
$\Ti(K_\chi)=H^{i}(\bo^0,K_{{\rm d}\chi})$ in $Vec_K$. The
algebras $\bo^0$ resp. $\tor^0$ resp. $\un^0$ are direct products
of scalar extensions of $\bo$ resp. $\tor$ resp. $\un$.  In
particular, $\tor^0\subseteq\bo^0$ is a reductive subalgebra
whence [HS], Thm. 12 implies that there is an isomorphism
\[H^{i}(\bo^0, K_{{\rm d}\chi})\simeq\sum_{j+k=i} H^{j}(\tor^0,K)\otimes_K
H^k(\bo^0,\tor^0,K_{{\rm d}\chi})\] in $Vec_K$ where
$H^{*}(\bo^0,\tor^0,K_{{\rm d}\chi})$ is the relative Lie algebra
cohomology with respect to $\tor^0$. Now $\tor^0$ is even toral
whence the argument preceding [loc.cit.], Thm. 13 implies that
$H^{*}(\bo^0,\tor^0,K_{{\rm d}\chi})\simeq H^*(\un^0,K_{{\rm
d}\chi})^{\tor^0}$. Finally, since $\tor^0$ is abelian, the
differential in ${\rm Hom}_K(\dot{\bigwedge}\tor^0,K)$ vanishes
identically and so one obtains $H^{j}(\tor^0,K)={\rm
Hom}_K(\bigwedge^{j}\tor^0,K)$.
\end{pr}
Remark: One has $H^k(\un^0,K_{{\rm d}\chi})=H^k(\un^0,K)$ in
$Vec_K$ and since $\un^0$ is nilpotent [D], Thm. 2 implies that
$\dim_K H^k(\un^0,K)\geq 1$ for $0\leq k\leq {\dim}_K\un ^0$.

\begin{cor}\label{Emerton}
If $\chi$ is smooth one has $R^1F_{\Q_p}^L({\rm
Ind}_{P_0}^{G_0}K_\chi)\neq 0$.
\end{cor}
\begin{pr}
We have ${\rm d}\chi=0$ whence $H^0(\un^0,K_{{\rm
d}\chi})^{\tor^0}=K$. By the corollary there is an injection
$\bigwedge^{1}\tor^0\rightarrow R^1F_{\Q_p}^L(K_\chi)$ whence the
claim follows from Cor. \ref{vertausch}.
\end{pr}

%%%%%%%%%%%%%%%%%%%%%%%%%%%%%%%%%%%%%%%%%%%%%%%%%%%%%%%%%%%%%%%%%

Tobias Schmidt

Mathematisches Institut

Westf\"alische Wilhelms-Universit\"at M\"unster

Einsteinstr. 62

D-48149 M\"unster, Germany

toschmid@math.uni-muenster.de

\end{document}